%% file: 0511uvSystemArx.tex
\documentclass[reqno,a4paper]{amsart}
\usepackage{amssymb}
\usepackage{latexsym}
\usepackage{amsmath}
\usepackage{amsfonts}
\usepackage[dvipdfmx]{hyperref}
\usepackage[dvipdfmx]{graphicx}

\numberwithin{equation}{section}

\newcommand{\pa}{\partial}
\newcommand{\ph}{\varphi}
\newcommand{\ep}{\varepsilon}
\newcommand{\R}{{\mathbb R}}

\newcommand{\supp}{\mathrm{supp}\,}
\newcommand{\lb}{\langle}
\newcommand{\rb}{\rangle}

\newtheorem{theorem}{Theorem}[section]

\newtheorem{lemma}[theorem]{Lemma}
\newtheorem{proposition}[theorem]{Proposition}
\theoremstyle{definition}

\theoremstyle{remark}
\newtheorem{remark}[theorem]{Remark}

\title[$p$-$q$ system]{Life span of small solutions to a system of wave equations}
\author[K.\,Hidano]{Kunio Hidano}
\author[K.\,Yokoyama]{Kazuyoshi Yokoyama}
\address{Department of Mathematics, Faculty of Education, 
Mie University, 1577 Kurima-machiya-cho Tsu, Mie Prefecture 514-8507, Japan}
\email{hidano@edu.mie-u.ac.jp}
\address{Department of Electrical and Electronic Engineering,
Faculty of Engineering\\
Hokkaido University of Science, 
7-15-4-1 Maeda, Teine, Sapporo, Hokkaido 006-8585, Japan}
\email{yokoyama@hus.ac.jp}


\begin{document}

\renewcommand{\thefootnote}{\fnsymbol{footnote}}
\footnote[0]{2010\textit{ Mathematics Subject Classification}.
 Primary 35L05, 35L15; Secondary 35L71}

\baselineskip=0.55cm
\begin{abstract}
We study the Cauchy problem with small initial data for a system of semilinear wave equations
%of the form
$\square u = |v|^p$, $\square v = |\pa_t u|^p$ in $n$-dimensional space. 
When $n \geq 2$, we prove that blow-up can occur for arbitrarily small data
if $(p, q)$ lies below a curve in $p$-$q$ plane. 
On the other hand, we show a global existence result for $n=3$
which asserts that a portion of the curve is in fact the borderline between global-in-time existence
and finite time blow-up. 
We also estimate the maximal existence time and get an upper bound, 
which is sharp at least for $(n, p, q)=(2, 2, 2)$
and $(3, 2, 2)$.
\end{abstract}

\maketitle

\input{intro.tex}

\input{blowup.tex}

\input{global.tex}

\input{uvSystem.tex}

\section*{Acknowledgements}
The authors wish to thank Chengbo Wang for 
his suggestion in the early stages of this research.
The first author was supported by Mie University Research Grant 2014 (Step-up (C)).

%\section*{References}

\bibliographystyle{amsplain}
%\bibliography{reference}
\input{0511uvSystemArx.bbl}

\end{document}

%% file: intro.tex
\section{Introduction}

We are interested in the study of systems of semilinear wave equations of the form
\begin{equation}
\label{b1}
 \square u = |v|^q,\quad
 \square v = |\pa_t u|^p.
\end{equation}
Here the unknown real-valued functions $u$, $v$ depend on $(t, x) \in 
[0,T) \times \R^n$ for some $T \in (0, \infty]$. 
Throughout this paper, we suppose $n \geq 2$. 
Given some functions $f$, $g$, $\tilde{f}$, $\tilde{g}$,
we consider the Cauchy problem for (\ref{b1}) with initial data
\begin{equation}
\label{b2}
 u(0,x)=\ep f(x),\quad \pa_t u(0,x)=\ep g(x),\quad
 v(0,x)=\ep \tilde{f}(x),\quad \pa_t v(0,x)=\ep \tilde{g}(x),
\end{equation}
where $\ep > 0$ is small. 
We want to determine, in the set $(1, \infty)^2$ 
of index-pairs $(p, q)$, the borderline
between global-in-time existence and finite time blow-up
for (\ref{b1})--(\ref{b2}) when $\ep$ is small. 

The system (\ref{b1}) reminds us of some related problems. 
Among them, the Cauchy problem for
\begin{equation}
\label{b10}
 \square u = |u|^q
\end{equation}
has been especially well-studied since John's pioneering work \cite{J79}. If $q > q_0(n)$ where
\[
 q_0(n) = \frac{n+1+\sqrt{n^2+10n-7}}{2(n-1)},
\]
then the Cauchy problem for (\ref{b10}) has a unique global-in-time 
solution for small and smooth initial data with compact support. 
If $1 < q \leq q_0(n)$ on the contrary, then there exists
$(f, g) \in C_0^{\infty}(\R^n)^2$ such that the Cauchy problem for (\ref{b10})
with initia data $(u(0), \pa_tu(0)) =
(\ep f, \ep g)$ has a solution which blows up
in finite time no matter how small $\ep > 0$ is. 
This phenomenon has been verified for all $n \geq 2$ through efforts
of several authors: see, e.g., 
\cite{J79}, \cite{Gl2}, \cite{Z95}, %\cite{KK95}, 
\cite{KK98},
\cite{LS}, \cite{GLS} for global existence results, 
\cite{J79}, \cite{Gl1}, \cite{S84}, \cite{Sc85}, \cite{R}, \cite{YZ06}
for blow-up results. Thus we have already understood that 
$q_0(n)$ is the critical exponent which separates 
global-in-time solvability and unsolvability for (\ref{b10}).

The similar phenomenon can be observed for the equation
\begin{equation}
\label{b11}
 \square v = |\pa_t v|^p.
\end{equation}
For this equation, the critical value of $p$ is
\[
 p_0(n) = \frac{n+1}{n-1}.
\]
In the case $1 < p \leq p_0(n)$, we know that solutions blow up
in finite time in general, no matter how small the initial data are
chosen. We also know that (\ref{b11}) with $p > p_0(n)$ has a unique global-in-time
solution for given small and smooth data, though the radial symmetry 
assumption is needed for $n \geq 4$ so far.
See, e.g., \cite{HT}, \cite{Tz}, %\cite{S83}, 
\cite{HWY2014} for global existence results,
\cite{A}, \cite{J81}, %\cite{Sc86}, \cite{R}, 
\cite{Z01} for blow-up results. 

Recently, combined effects of the two nonlinearities above have been studied
by \cite{HZ} and \cite{HWY2014}. They considered the Cauchy problem for 
\begin{equation}
\label{b12}
 \square u = |u|^q + |\pa_t u|^p
\end{equation}
(only when $n=2, 3$ in \cite{HWY2014}).
According to their results, we can conclude that
\begin{equation}
\label{b14}
 q = \frac{4}{(n-1)p-2}+1 \quad (q>q_0(n),\ p>p_0(n))
\end{equation}
is a portion of the critical curve in the $p$-$q$ plane of index-pairs $(p, q)$
(see \cite{HZ}, \cite{HWY2014} for details).
Observe that $(p_0(n), q_0(n))$ belongs to the domain
$q < 4/\{(n-1)p-2\}+1$, in which $(p, q)$
corresponds to the blow-up case.
As is pointed out in \cite{HZ}, this means that
there exists $(p, q)$ such that a solution of the Cauchy problem 
(\ref{b12}) blows up in finite time for arbitrarily small initial data, 
while the Cauchy problems for (\ref{b1}) and (\ref{b11}) with
small data have global solutions.
Another remarkable feature of (\ref{b12}) is that,
unlike (\ref{b1}) and (\ref{b11}), 
we can establish global existence for $(p, q)$ being on the
curve (\ref{b14}) (see \cite{HWY2014}).

When it comes to systems of wave equations, 
the Cauchy problem for
\begin{equation}
 \square u = |v|^q,\quad \square v = |u|^p
\label{b15}
\end{equation}
has been well-studied. It is known that 
\begin{equation}
\label{b13}
 \max \left\{\frac{p+2+q^{-1}}{pq-1},\ 
\frac{q+2+p^{-1}}{pq-1}\right\}-\frac{n-1}{2}=0
\end{equation}
is the critical curve in the $p$-$q$ plane of index-pairs $(p, q)$.
See, e.g., 
   %\cite{AKT}, %\cite{D98}, 
   \cite{DGM}, \cite{DM}, 
   \cite{KO99}, \cite{GTZ}. %\cite{KO00}, 
   %\cite{Ku05},
   %\cite{KTW}, 
Remark that $(q_0(n), q_0(n))$ is on the critical curve. Also,
the curve (\ref{b13}) lies above the curve (\ref{b14})
when $q > q_0(n)$.

Now let us turn to the problem (\ref{b1})--(\ref{b2}). 
We first show a blow-up result which is valid for all $n \geq 2$. 
When discussing blow-up, 
we say that $(u, v)$ satisfies (\ref{b1}) if the following two conditions
(i), (ii) hold:
\begin{enumerate}
\renewcommand{\labelenumi}{(\roman{enumi})}
\item 
The equalities
\begin{align}
\label{y1}
 &\int_0^T \lb |v(s,\cdot)|^q,\psi(s,\cdot)\rb\,ds
=\int_0^T \lb u(s,\cdot), (\pa_s^2 \psi -\Delta \psi)(s,\cdot)\rb\,ds,\\
 &\int_0^T \lb |\pa_s u(s,\cdot)|^p,\psi(s,\cdot)\rb\,ds
=\int_0^T \lb v(s,\cdot), (\pa_s^2 \psi -\Delta \psi)(s,\cdot)\rb\,ds
\end{align}
hold for any $\psi \in C_0^{\infty}((0,T)\times \R^n)$. 
\item 
$u$ and $v$ have the continuity and integrability properties
\begin{align}
\label{y2}
 & u, v, \pa_t u, \pa_t v \in C([0,T);\,L^1(\R^n)),\\
\label{y3}
 & v \in C((0,T);\,L^q(\R^n)),\quad
 \pa_t u \in C((0,T);\,L^p(\R^n)).
\end{align}
\end{enumerate}
In particular, the initial conditions (\ref{b2}) should be satisfied
in the sense of (\ref{y2}).
\begin{theorem}\label{y4}
 Let $n \geq 2$, $1<q$, $1 < p < 2n/(n-1)$ and
\begin{equation}
\label{y5}
 \left(\frac{n-1}{2}p-1\right)(pq-1) < p+2.
\end{equation}
Assume that the initial value problem {\rm (\ref{b1})--(\ref{b2})}
with non-negative initial data has a solution
$(u, v)$ in the sense described above, which is localized inside a cone:
\begin{equation}
\label{y6}
 \supp (u, v, \pa_t u, \pa_t v) \subset \{(t,\,x) \in [0, T) \times \R^n\,:\,
|x| \leq t+1\}.
\end{equation}
We also assume that $f+g$ and $\tilde{g}$ do not vanish identically. 
Then there exists a positive constant $C_1$,
depending only on $n, p, q, f, g, \tilde{f}, \tilde{g}$, such that
\begin{equation}
\label{y7}
 T \leq C_1 \ep^{-\frac{p(pq-1)}{p+2-\left(\frac{n-1}{2}p-1\right)(pq-1)}}.
\end{equation}
\end{theorem}

Observe that $(p_0(n), q_0(n))$ belongs to the domain
(\ref{y5}). Thus the same remark as the one stated on (\ref{b12}) above
applies to (\ref{b1}).
We infer from Theorem \ref{y4} that a portion of the curve
\begin{equation}
\label{b5}
 \left(\frac{n-1}{2}p-1\right)(pq-1) = p+2
\end{equation}
gives a component part of the critical curve in the $p$-$q$ plane.
In order to confirm this expectation,
we should establish a global existence theorem for $(p, q)$
above the curve (\ref{b5}). 
We solve the system of the integral equations
\begin{align}
\label{z1}
 u(t,r) & = \ep U_0(t,r)+ L|v|^q (t,r), \\
 \label{z2}
 v(t,r) & = \ep V_0(t,r)+ L|\pa_t u|^p (t,r),
\end{align}
where
\begin{align}
 U_0(t,r)&=\frac{1}{2r}\left\{(r+t)f(r+t)+(r-t)f(r-t)\right\}
+\frac{1}{2r}\int_{r-t}^{r+t}\rho g(\rho)\,d\rho,
\nonumber\\
 V_0(t,r)&=\frac{1}{2r}\left\{(r+t)\tilde{f}(r+t)+(r-t)\tilde{f}(r-t)\right\}
+\frac{1}{2r}\int_{r-t}^{r+t}\rho \tilde{g}(\rho)\,d\rho
\nonumber
\end{align}
and
\begin{equation}
 LF(t,r) = \frac{1}{2r}\int_0^t \int_{r-(t-s)}^{r+(t-s)}\rho 
  F(s,\rho)\,d\rho ds.
\nonumber
\end{equation}
The system (\ref{z1})--(\ref{z2}) is a radially symmetric version of the 
original system with $n=3$,
and we naturally suppose that $f(r)$, $g(r)$, $u(t,r)$, $\ldots$ are even functions of $r$. 
Hence, the lower limits of the integrals above may be replaced by $|r-t|$ or $|r-(t-s)|$. 
$LF(t,r)$ is also an even function of $r$. 

Throughout the proof of global existence, we suppose that
\begin{equation}
\label{z3}
 q > 2, \quad 2 < p < 3.
\end{equation}
We want to show that if in addition
\begin{equation}
\label{z4}
 (p-1)(pq-1)>p+2
\end{equation}
is assumed, then the system (\ref{z1})--(\ref{z2}) has
a global solution for sufficiently small data.

For this purpose, we introduce a weighted norm.
We define weight functions $w_1$, $w_2$, $w_3$ by
\begin{align}
 w_1(t,r)&=\lb r\rb \lb t-r\rb^{\mu/p}, \\
 \label{z5}
 w_2(t,r)&=
 \begin{cases}
  \lb r \rb^{p-2} \lb t+r\rb^{\mu} & \left(r < t/2\right) \\
  \lb t-r\rb^{p-3+\mu}\lb t+r\rb & \left(r \geq t/2\right)
 \end{cases},\\
 w_3(t,r)&= \lb t-r\rb^{\mu}
\end{align}
for $t \geq 0$ and $r \geq 0$. Here we use the concise notation $\lb \xi\rb=\sqrt{1+|\xi|^2}$ as usual, but
we may replace it by $1+|\xi|$ below.
By (\ref{z4}), we can choose $\mu$ so that 
\begin{equation}
\label{z6}
 3-p < \mu < 1,\quad 
p\left\{2-(p-2)q\right\} < (pq-1)\mu
\end{equation}
are satisfied. Then we define
\begin{equation}
 \|(w, v)\|
=\|w_1w\|_{L^{\infty}_{t,r}}
+\|w_2v\|_{L^{\infty}_{t,r}}
+\|w_3r\pa_r v\|_{L^{\infty}_{t,r}}
\nonumber
\end{equation}
where $L^{\infty}_{t,r} = L^{\infty}([0, \infty) \times [0, \infty))$.

\begin{theorem}\label{z7}
Suppose {\rm (\ref{z3})}, {\rm (\ref{z4})}. 
Suppose also that $f, \tilde{f} \in C_0^2(\R)$,\ $g, \tilde{g} \in C_0^1(\R)$ are even functions. 
There exist positive numbers $\ep_2$, $C_2$ such that
we have a unique global solution $(u, v)$ of {\rm (\ref{z1})--(\ref{z2})}
for $0 < \ep < \ep_2$,
satisfying $\|(\pa_tu, v)\| \leq C_2\ep$ and
\begin{align}
 &u(t,-r)=u(t,r),\ v(t,-r)=v(t,r),\nonumber\\
 &(u, v) \in C([0, \infty)\times \R) \times C([0,\infty)\times \R), 
\quad \pa_t u \in C([0, \infty)\times \R), \nonumber\\
  & (ru, rv) \in C^2([0,\infty)\times \R) \times C^2([0,\infty)\times \R). 
\nonumber
\end{align}
\end{theorem}
Thus we may say that a portion of the curve (\ref{b5}) is truely the  borderline 
between global-in-time existence and finite time blow-up
for (\ref{b1})--(\ref{b2}). 
It is known for (\ref{b15})
that finite time blow-up occurs for $(p, q)$ being on the critical curve,
but this problem is poorly understood for (\ref{b1}) and
we have no results about the case where $(p, q)$ belongs to 
 the critical curve at this time. 
Note that the curve (\ref{b5}) lies below the curve (\ref{b14}) for $p > p_0(n)$. 

If radial symmetry is not assumed,
we can use the method of \cite{H98}, \cite{HWY2014} to get 
some partial results of global existence. However, the present
authors have not yet suceeded in
obtaining a global existence result right up to the critical curve. 
Instead, we treat the case $(p, q) = (2, 2)$ for $n = 2, 3$
to describe the method of \cite{HWY2014}. It is of some interest
since we can solve (\ref{b1})--(\ref{b2}) up to 
$T=c\ep^{-\frac{6}{10-3n}}$ in this way (recall (\ref{y7})).

We assume that
\begin{gather}
\label{b8}
 f \in L^2(\R^n),\quad x g \in L^2 (\R^n), \quad 
x \tilde{f} \in \dot{H}^{1/4}(\R^n),\\
 x^a \pa_x^b \pa_x f \in L^2 (\R^n), \quad 
 x^a \pa_x^b g \in L^2 (\R^n),\\
 x^a \pa_x^b \tilde{f} \in \dot{H}^{1/4}(\R^n),\quad
 x^a \pa_x^b \tilde{g} \in \dot{H}^{-3/4}(\R^n)
\label{b9}
\end{gather}
for $|a| \leq |b| \leq 2$, where
${\dot H}^s({\mathbb R}^n)$
is the standard homogeneous Sobolev space. We set
\begin{align}
\label{b16}
\Lambda_1&
:=\|f\|_{L^2}+\sum_{i=1}^n\|x_ig\|_{L^2}
+\sum_{|a|\leq |b| \leq 2}\left(
\|x^a\pa_x^b \pa_x f\|_{L^2}
+\|x^a\pa_x^b g\|_{L^2}\right)\\
\intertext{and}
\label{b17}
\Lambda_2&
:=\sum_{i=1}^n\|x_i\tilde{f}\|_{\dot{H}^{1/4}}
+\sum_{|a|\leq |b|\leq 2}
\left(
\|x^a\pa_x^b{\tilde f}\|_{\dot{H}^{1/4}}
+\|x^a\pa_x^b{\tilde g}\|_{\dot{H}^{-3/4}}
\right).
\end{align}

Let us denote $\pa_0 = \pa_t$ as usual. Following Klainerman \cite{Kl85}--\cite{Kl87}, 
we introduce several partial differential operators as follows: 
$L_j=t\partial_j+x_j\partial_0 $ $(j=1,\dots,n)$, 
$\Omega_{kl}=x_k\partial_l-x_l\partial_k$ $(1\leq k<l\leq n)$, 
$L_0=t\partial_0+x_1\partial_1+\cdots+x_n\partial_n$. 
Operators 
$\partial_0,\dots,\partial_n$, $L_1,\dots,L_n$, $\Omega_{12},\dots,
\Omega_{1n},\Omega_{23},\dots,\Omega_{n-1n}$ and $L_0$ are denoted by 
$\Gamma_0,\dots,\Gamma_\nu$ in this order, where 
$\nu:=(n^2+3n+2)/2$. 
For a multi-index $\alpha=(\alpha_0,\dots,\alpha_\nu)$, 
$\Gamma^{\alpha_0}_0\cdots\Gamma^{\alpha_\nu}_\nu$ is denoted by 
$\Gamma^\alpha$.

\begin{theorem}\label{b6}
 Let $(p, q)=(2, 2)$ and $n=2, 3$. Assume {\rm (\ref{b8})--(\ref{b9})}.
Then, there exist $\ep_3$, $C_3$ depending on $n$, $\Lambda_1$ and
$\Lambda_2$ such that the Cauchy problem {\rm (\ref{b1})--(\ref{b2})}
with $0 < \ep < \ep_3$ and $T = C_3\ep^{-\frac{6}{10-3n}}$
admits a unique solution satisfying
\begin{eqnarray}
& &
\Gamma^\alpha u\in C([0,T];H^1({\mathbb R}^n)),\quad
\partial_t\Gamma^\alpha u\in C([0,T];L^2({\mathbb R}^n)),
\nonumber\\
& &
\Gamma^\alpha v\in C([0,T];\dot{H}^{1/4}({\mathbb R}^n)),
\nonumber
\end{eqnarray}
for $|\alpha| \leq 2$ and
\begin{eqnarray*}
& &
\sup_{0<t<T}
(1+t)^{-1}\|u(t)\|_{L^2}
+
\sum_{|\alpha|\leq 2\atop 0\leq j\leq n}
\sup_{0<t<T}
\|\partial_j\Gamma^\alpha u(t)\|_{L^2}\\
& &
\hspace{2cm}
+
\sum_{|\alpha|\leq 2}
\sup_{0<t<T}
(1+t)^{-1/12}\|\Gamma^\alpha v(t)\|_{\dot{H}^{1/4}}
\leq C\ep.\nonumber
\end{eqnarray*}
\end{theorem}

In the rest of this paper, we prove Theorem \ref{y4}, 
Theorem \ref{z7} and Theorem \ref{b6}, in Sections 2, 3 and 4, respectively.

%% file: blowup.tex
\section{Proof of Theorem \ref{y4}}

Our aim in this section is to prove Theorem \ref{y4}.
Let us define
\begin{equation}
 F(t)=\int_{\R^n}u(t,x)\,dx,\quad
 G(t)=\int_{\R^n}v(t,x)\,dx.
\end{equation}
We adopt the strategy of deriving a system of ordinary differential inequalities
with respect to $F(t)$ and $G(t)$,
which causes blow-up of solutions. 
Remark that if $(u, v)$ solves the Cauchy problem in the sense
stated above Theorem \ref{y4}, 
then we can see that $F''(t)=\|v(t,\cdot)\|_{L^q}^q$ and 
$G''(t)=\|\pa_tu(t,\cdot)\|_{L^p}^p$ for $t \in (0, T)$. 
Hence we conclude from (\ref{y2})--(\ref{y3}) that $F, G \in C^2((0,T)) 
\cap C^1([0, T))$. Furthermore, $F'(t)$ and $G(t)$ are non-negative.
Thus the following lemma is crucial for the proof of Theorem \ref{y4}.
\begin{lemma}\label{y8}
Let $F, G \in C^2((0, T)) \cap C^1([0, T))$. We suppose that
they satisfy
$F'(t), G(t) \geq 0$ and
\begin{align}
\label{y9}
 &F''(t) \geq A(t+1)^{-\alpha}G(t)^q,\\
\label{y10}
 &G''(t) \geq B(t+1)^{-\beta}F'(t)^p
\end{align}
for all $t \in (0,T)$,
with some positive constants $A, B$ and exponents
$p, q > 1$, $\alpha, \beta \geq 0$.
Assume also $G'(0) > 0$ and
\begin{align}
\label{y11}
 &G(t) \geq \kappa t^a
\end{align}
on some sub-interval $[t_0, T)$, where $\kappa, a > 0$ are constants. 
If $\max \left\{t_0, G(0)/G'(0), 1\right\}$ $< T/2$ and
\begin{equation}
 \beta + \alpha p < p+2 + a(pq-1),
\label{y12}
\end{equation}
then 
\begin{equation}
 T \leq C\kappa^{-\frac{pq-1}{p+2+a(pq-1)-(\beta + \alpha p)}},
\end{equation}
where $C$ is a constant depending only on $A, B, p, q, \alpha, \beta, a$.
\end{lemma}

\subsection{Proof of Lemma {\rm \ref{y8}}}
We follow the method of Kurokawa, Takamura and Wakasa (see 
Lemma 2.1 of \cite{KTW}).
Note beforehand that we have $G'(t) > 0$ as a result of $G'(0)>0$ and (\ref{y10}).
Hence by (\ref{y9}), we have
\begin{align}
 \int_0^t F''(s)G'(s)\,ds 
 &\geq A\int_0^t (s+1)^{-\alpha}G(s)^q G'(s)\,ds
\nonumber\\
 &\geq A(t+1)^{-\alpha}\int_0^t G(s)^q G'(s)\,ds.
\nonumber
\end{align}
Applying the integration by parts formula on the left-hand side, 
we get
\begin{equation}
 F'(t)G'(t) \geq \frac{A}{q+1}(t+1)^{-\alpha}\left(
G(t)^{q+1}-G(0)^{q+1}\right)
\label{y13}
\end{equation}
for $t > 0$.

We turn to the inequality (\ref{y10}) next. 
Multiplying both sides of it by $G'(t)^p$ and using (\ref{y13}), 
we get
\begin{align}
 G'(t)^pG''(t) &\geq B(t+1)^{-\beta}\left(
\frac{A}{q+1}(t+1)^{-\alpha}\left(
G(t)^{q+1}-G(0)^{q+1}\right)
\right)^p \nonumber
\end{align}
for $t > 0$. Thus we immediately obtain
\begin{align}
 \int_0^tG'(s)^pG''(s)\,ds 
&\geq \frac{A^pB}{(q+1)^p}(t+1)^{-\beta-\alpha p}\int_0^t
\left(G(s)^{q+1}-G(0)^{q+1}\right)^p\,ds.
\nonumber
\end{align}
We can directly compute the left-hand side of this inequality. 
As for the right-hand side, we estimate the integral as follows:
\begin{align}
 \lefteqn{\int_0^t \left(G(s)^{q+1}-G(0)^{q+1}\right)^p\,ds}
\nonumber\\
 &= \int_0^t \left(G(s)^{q+1}-G(0)^{q+1}\right)^p
\frac{\left(G(s)^{q+1}\right)'}{(q+1)G(s)^q G'(s)}\,ds
\nonumber\\
 &\geq \frac{1}{(q+1)G(t)^q G'(t)}\int_0^t \left(G(s)^{q+1}-G(0)^{q+1}\right)^p
\left(G(s)^{q+1}\right)'\,ds.
\nonumber
\end{align}
In the last inequality, we have used the fact that $G(s)$ and $G'(s)$ are 
increasing functions. 
Combining the estimates above, we have
\begin{equation}
 G'(t)^{p+2}
\geq \frac{A^pB}{(q+1)^{p+1}}(t+1)^{-\beta-\alpha p}
\frac{\left(G(t)^{q+1}-G(0)^{q+1}\right)^{p+1}}{G(t)^q}
\label{y14}
\end{equation}
for $t > 0$.

Now we suppose $t \geq G(0)/G'(0)$. 
Then we see $G(t) \geq 2G(0)$, 
because
\begin{equation}
 G(t) = G(0) + \int_0^t G'(s)\,ds \geq G(0) + tG'(0)
\geq 2G(0).
\nonumber
\end{equation}
Hence (\ref{y14}) implies
\begin{align}
 G'(t)^{p+2} &
\geq \frac{A^pB}{(q+1)^{p+1}}(t+1)^{-\beta-\alpha p}
(1-2^{-q-1})^{p+1}G(t)^{(q+1)(p+1)-q}
\nonumber
\end{align}
for $t \geq G(0)/G'(0)$. 
Since $t+1 \leq 2t$ if $t \geq 1$, we further get
\begin{equation}
 G'(t)
\geq Dt^{-\frac{\beta+\alpha p}{p+2}}
G(t)^{\frac{pq+p+1}{p+2}}
\label{y15}
\end{equation}
for $t \geq \max \left\{G(0)/G'(0), 1\right\}$, 
where 
\begin{equation}
 D=\left\{\frac{A^pB(1-2^{-q-1})^{p+1}}{(q+1)^{p+1}2^{\beta + 
    \alpha p}}\right\}^{\frac{1}{p+2}}.
\nonumber
\end{equation}

\bigskip

In what follows, we will show that the inequality (\ref{y15}) 
cannot hold if $t$ is large enough. 
We use (\ref{y11}) in this position. 
By (\ref{y15}) and (\ref{y11}), 
we have
\begin{align}
 G(t)^{-1-\delta}G'(t) &\geq Dt^{-\frac{\beta+\alpha p}{p+2}}G(t)^{\frac{pq-1}{p+2}-\delta}
\nonumber \\
 &\geq D\kappa^{\frac{pq-1}{p+2}-\delta}t^{\frac{a(pq-1)-(\beta+\alpha p)}{p+2}-a\delta}
\label{y16}
\end{align}
for $t \geq \max \{t_0, G(0)/G'(0), 1\}$,
where we choose $\delta$ so that (recall (\ref{y12}))
\[
 0 < \delta \leq \frac{pq-1}{p+2}, \quad \frac{a(pq-1)-(\beta+\alpha 
 p)}{p+2}-a\delta > -1.
\]
Now we take arbitrary $\tilde{T} < T$ so that
$\max \{t_0, G(0)/G'(0), 1\} \leq \tilde{T}/2$. 
Integrating (\ref{y16}) over $\left[\tilde{T}/2, \tilde{T}\right]$,
we obtain
\begin{equation}
 \frac{1}{\delta}G\left(\frac{\tilde{T}}{2}\right)^{-\delta} 
\geq \frac{D\kappa^{\frac{pq-1}{p+2}-\delta}%
(1-2^{-1-\frac{a(pq-1)-(\beta+\alpha p)}{p+2}+a\delta})}{%
1+\frac{a(pq-1)-(\beta+\alpha p)}{p+2}-a\delta}\tilde{T}^{1+\frac{a(pq-1)-(\beta+\alpha p)}{p+2}-a\delta}.
\label{y17}
\end{equation}
Using (\ref{y11}) again, we see that
\begin{equation}
 \mbox{LHS of (\ref{y17})} \leq \frac{1}{\delta}\kappa^{-\delta}\left(\frac{\tilde{T}}{2}\right)^{-a\delta}.
\label{y18}
\end{equation}
Hence, it follows from (\ref{y17}) and (\ref{y18}) that
\begin{equation}
 \tilde{T} \leq C \kappa^{-\frac{pq-1}{p+2+a(pq-1)-(\beta + \alpha p)}},
\end{equation}
where $C$ is a constant depending only on $A, B, p, q, \alpha, \beta, a$.
This completes the proof of Lemma \ref{y8}.
\hfill $\square$

\subsection{Proof of Theorem {\rm \ref{y4}}}

Following Yordanov--Zhang \cite{YZ05}, we introduce the two positive functions
\begin{equation}
 \ph_1(x) = \int_{S^{n-1}}e^{x\cdot \omega}\,d\omega,\quad
 \psi_1(t,x) = \ph_1(x)e^{-t}.
\nonumber
\end{equation}
Note that $\psi_1$ satisfies $\square \psi_1 = 0$. 
By a standard manner, we can deduce from (\ref{y1})--(\ref{y3})
that
\begin{align}
\nonumber
 \int_0^t \lb |v(s,\cdot)|^q, \psi(s,\cdot)\rb\,ds
&= \lb \pa_t u(s,\cdot), \psi(s,\cdot)\rb\Big|_{s=0}^{s=t}
-\lb u(s,\cdot), \pa_t\psi(s,\cdot)\rb\Big|_{s=0}^{s=t}\\
&\quad +\int_0^t \lb u(s,\cdot),(\pa_t^2\psi -\Delta \psi)(s,\cdot)\rb\,ds
\nonumber
\end{align}
for any $\psi \in C^{\infty}([0, T) \times \R^n)$ such that $\psi (s,\cdot) 
\in C_0^{\infty}(\R^n)$ for each $s \in [0, T)$. Thanks to (\ref{y6}),
we may substitute $\psi = \psi_1$. As a result, we easily obtain
\begin{equation}
 \int_0^t \lb |v(s,\cdot)|^q, \psi_1(s,\cdot)\rb\,ds
= \lb \pa_t u(t,\cdot) + u(t,\cdot), \psi_1(t,\cdot)\rb
-\ep \lb f + g, \ph_1\rb.
\nonumber
\end{equation}
Moreover, we see from (\ref{y1}) that
\begin{align}
\nonumber
 \lefteqn{\frac{d}{dt}\lb \pa_t u(t,\cdot), \psi_1(t,\cdot)\rb}\\
 & = \lb |v(t,\cdot)|^q, \psi_1(t,\cdot)\rb + \lb u(t,\cdot),\Delta \psi_1(t,\cdot)\rb
-\lb \pa_t u(t,\cdot),\psi_1(t,\cdot)\rb \nonumber\\
 & = \lb |v(t,\cdot)|^q, \psi_1(t,\cdot)\rb + \lb u(t,\cdot), \psi_1(t,\cdot)\rb
-\lb \pa_t u(t,\cdot),\psi_1(t,\cdot)\rb \nonumber
\end{align}
for $t \in (0, T)$. Combining these estimates, we get
\begin{align}
e^{-2t}\frac{d}{dt}\left(e^{2t}\lb \pa_t u(t,\cdot), \psi_1(t,\cdot)\rb\right)
&= \frac{d}{dt}\lb \pa_t u(t,\cdot), \psi_1(t,\cdot)\rb
+ 2\lb \pa_t u(t,\cdot),\psi_1(t,\cdot)\rb\nonumber\\
  &\geq \ep \lb f+g, \ph_1\rb,\nonumber
\end{align}
which gives
\begin{equation}
\label{y19}
 \lb \pa_t u(t,\cdot), \psi_1(t,\cdot)\rb \geq \frac{1}{4}\ep \lb f+g, \ph_1\rb
\end{equation}
for $t \geq \log \sqrt{2}$.

On the other hand, according to Lemma 2.3 of Yordanov--Zhang \cite{YZ05}, 
\begin{equation}
\label{y20}
 \|\psi_1(t,\cdot) \|_{L^{p'}(|x| \leq t+1)} \leq C(t+1)^{(n-1)\left(\frac{1}{2}-\frac{1}{p}\right)}
\end{equation}
for $t \geq 0$ (see also Zhou--Han \cite{ZH14}). Hence, by (\ref{y19}) and (\ref{y20}), 
we have
\begin{equation}
\nonumber
 G''(t) = \|\pa_t u(t,\cdot)\|_{L^p}^p \geq C\ep^p(t+1)^{-\frac{(n-1)(p-2)}{2}}
\end{equation}
for $t \geq \log \sqrt{2}$. Since $(n-1)(p-2)/2<1$ by assumption, 
we eventually get the estimate
\begin{equation}
\nonumber
 G(t) \geq C\ep^p (t+1)^{2-\frac{(n-1)(p-2)}{2}}
\end{equation}
for $t \geq t_0$, if $t_0$ is large enough.

Meanwhile, by (\ref{y6}) and the H\"older inequality, we 
immediately have
\begin{align}
 &\left(F'(t)\right)^p \leq C(t+1)^{n(p-1)}\|\pa_tu(t,\cdot)\|_{L^p}^p
= C(t+1)^{n(p-1)}G''(t),\nonumber\\
&\left(G(t)\right)^q \leq C(t+1)^{n(q-1)}\|v(t,\cdot)\|_{L^q}^q
= C(t+1)^{n(q-1)}F''(t)\nonumber
\end{align}
for $t > 0$.

Now we are ready to apply Lemma \ref{y8} to the present problem.
Set
\begin{equation}
\nonumber
 a = 2-\frac{(n-1)(p-2)}{2},\quad \alpha = n(q-1),\quad
\beta = n(p-1)
\end{equation}
and $\kappa = C\varepsilon^p$. 
Noting $\beta + \alpha p = n(pq-1)$, 
we see that (\ref{y12}) is equivalent to (\ref{y5}).
Hence we have (\ref{y7}) as a conclusion of Lemma \ref{y8}.
\qed

%% file: global.tex
\section{Proof of Theorem \ref{z7}}

To control the iteration procedure, we need the estimates of derivatives. 
Computing the first derivatives of (\ref{z1}) and (\ref{z2}), we find that
\begin{align}
 r\pa_t u(t, r) &= \ep r\pa_t u_0(t, r)+K_+|v|^q (t, r),
\nonumber\\
 \pa_r \left\{rv(t, r)\right\} &= \ep \pa_r \left\{rv_0(t, r)\right\}
+K_-|\pa_tu|^p (t, r),\nonumber
\end{align}
where
\begin{equation}
 K_{\pm}F(t, r) 
= \frac{1}{2}\int_0^t \left\{(r+t-s) F(s,r+t-s) \pm (r-t+s) 
F(s,r-t+s)\right\}\,ds. \nonumber
\end{equation}
As functions of $r$, 
$K_+F(t,r)$ is odd and $K_-F(t,r)$ is even. 

In order to prove Theorem \ref{z7}, we use the following lemmas.
Let us recall the definitions of the operator $L$ and the weights $w_1$, 
$w_2$, $w_3$ for convenience:
\begin{equation}
 LF(t,r) = \frac{1}{2r}\int_0^t \int_{r-(t-s)}^{r+(t-s)}\rho 
  F(s,\rho)\,d\rho ds,
\nonumber
\end{equation}
\begin{align*}
 w_1(t,r)&=\lb r\rb \lb t-r\rb^{\mu/p}, \\
 w_2(t,r)&=
 \begin{cases}
  \lb r \rb^{p-2} \lb t+r\rb^{\mu} & \left(r < t/2\right) \\
  \lb t-r\rb^{p-3+\mu}\lb t+r\rb & \left(r \geq t/2\right)
 \end{cases},\\
 w_3(t,r)&= \lb t-r\rb^{\mu}.
\end{align*}
\begin{lemma}\label{z8}
Suppose {\rm (\ref{z3})}, {\rm (\ref{z4})}. 
Then we have
 \begin{align}
 \label{z9}
&\|w_2L|w|^{p}\|_{L^{\infty}_{t,r}}
 \leq  C\|w_1 w\|^{p}_{L^{\infty}_{t,r}},\\
%%%
\label{z10}
&\|\lb t-r \rb^{\mu/p}K_+|v|^{q}\|_{L^{\infty}_{t,r}}
 \leq  C\|w_2 v\|^{q}_{L^{\infty}_{t,r}},\\
%%%
\label{z11}
&\left\|\lb t \rb^{\mu/p} r^{-1}
K_+|v|^{q}\right\|_{L^{\infty}_{t,r}(r < 1)}
 \leq  C\|w_2 v\|^{q}_{L^{\infty}_{t,r}}
+C\|w_2 v\|^{q-1}_{L^{\infty}_{t,r}}
\|w_3 r\pa_rv\|_{L^{\infty}_{t,r}},\\
%%%
\label{z12}
&\|w_3K_-|w|^{p}\|_{L^{\infty}_{t,r}}
 \leq  C\|w_1 w\|^{p}_{L^{\infty}_{t,r}}.
 \end{align}
\end{lemma}

\begin{lemma}\label{z13}
Suppose {\rm (\ref{z3})}, {\rm (\ref{z4})}. 
Then we have
 \begin{align}
&\|w_2L[|w|^{p}-|\bar{w}|^p]\|_{L^{\infty}_{t,r}}
 \leq  C\left(\|w_1 w\|^{p-1}_{L^{\infty}_{t,r}}
+\|w_1 \bar{w}\|^{p-1}_{L^{\infty}_{t,r}}
\right)
 \|w_1 (w-\bar{w})\|_{L^{\infty}_{t,r}},\\
%%% 
\lefteqn{\|\lb t-r \rb^{\mu/p}K_+[|v|^{q}-|\bar{v}|^q]\|_{L^{\infty}_{t,r}}}
\nonumber\\
& \leq  C\left(\|w_2 v\|^{q-1}_{L^{\infty}_{t,r}}
+\|w_2 \bar{v}\|^{q-1}_{L^{\infty}_{t,r}}
\right)
\|w_2 (v-\bar{v})\|_{L^{\infty}_{t,r}},
\label{z26}\\
%%%
\lefteqn{\left\|\lb t \rb^{\mu/p}r^{-1}K_+[|v|^{q}-|\bar{v}|^q]\right\|_{L^{\infty}_{t,r}(r < 1)}}
\nonumber\\
& \leq  C\left(\|w_2 v\|^{q-1}_{L^{\infty}_{t,r}}
+\|w_2 \bar{v}\|^{q-1}_{L^{\infty}_{t,r}}
\right)
\|w_2 (v-\bar{v})\|_{L^{\infty}_{t,r}}\nonumber\\
&\quad +C\left(\|w_2 v\|^{q-2}_{L^{\infty}_{t,r}}
+\|w_2 \bar{v}\|^{q-2}_{L^{\infty}_{t,r}}
\right)\|w_3r\pa_r v\|_{L^{\infty}_{t,r}}
\|w_2 (v-\bar{v})\|_{L^{\infty}_{t,r}}
\nonumber\\
&\quad+ C\|w_2\bar{v}\|^{q-1}_{L^{\infty}_{t,r}}
\|w_3r(\pa_r v - \pa_r \bar{v})\|_{L^{\infty}_{t,r}},
\label{z27}\\
%%%
&\|w_3K_-[|w|^{p}-|\bar{w}|^p]\|_{L^{\infty}_{t,r}}
 \leq  C\left(\|w_1 w\|^{p-1}_{L^{\infty}_{t,r}}
+\|w_1 \bar{w}\|^{p-1}_{L^{\infty}_{t,r}}
\right)
 \|w_1 (w-\bar{w})\|_{L^{\infty}_{t,r}}.
 \end{align}
\end{lemma}

{\it Proof of Theorem {\rm \ref{z7}}}.\ We consider the system of integral equations
\begin{align}
\label{z14}
 w(t,r) & = \ep \pa_tu_0(t,r)+ r^{-1}K_+|v|^q (t,r), \\
 \label{z15}
 v(t,r) & = \ep v_0(t,r)+ L|w|^p (t,r)
\end{align}
in the complete metric space
\begin{align}
 X_{\ep}=&\{(w, v) \in C([0,\infty)\times \R) \times C([0,\infty)\times \R)\,:\,\nonumber\\
  &\quad w(t,-r)=u(t,r),\ v(t,r)=v(t,-r),\nonumber\\
  &\quad  r\pa_r v \in C([0,\infty)\times \R),\nonumber\\
&\quad \|(w,v)\| =\|w_1w\|_{L^{\infty}_{t,r}}
+\|w_2v\|_{L^{\infty}_{t,r}}
+\|w_3r\pa_r v\|_{L^{\infty}_{t,r}} \leq C_2\ep\}.\nonumber
\end{align}
Lemmas \ref{z8}, \ref{z13} imply that
(\ref{z14})--(\ref{z15}) can be solved by
applying the Banach fixed-point theorem for the mapping
\begin{equation}
 (w, v) \longmapsto (\ep \pa_tu_0+ r^{-1}K_+|v|^q,\ 
 \ep v_0+ L|w|^p),\nonumber
\end{equation}
if $C_2$ is suitably chosen and $\ep$ is sufficiently small.
Since the argument is standard, we omit the details.
Setting
\begin{equation}
 u(t,r) := \ep f(r) + \int_0^tw(s,r)\,ds, \nonumber
\end{equation}
we see that $(u, v)$ satisfies {\rm (\ref{z1})--(\ref{z2})}. 
Moreover, we can easily check the regularity property by using the 
integral equations (\ref{z1})--(\ref{z2}).
\qed

We will prove Lemma \ref{z8} in the following subsections.
We omit the proof of Lemma \ref{z13}, since it is similar to the proof 
of Lemma \ref{z8}.
To prove Lemma \ref{z13}, we adapt the proof of
Lemma \ref{z8} by noting
\begin{align}
 &\big||v|^q-|\bar{v}|^q\big|
\leq C\left(|v|^{q-1}+|\bar{v}|^{q-1}\right)|v-\bar{v}|
\nonumber\\
\intertext{for (\ref{z26}),}
 &\lefteqn{
\big|
\left(|v(s,r^+)|^q-|\bar{v}(s,r^+)|^q\right)
-\left(|v(s,|r^-|)|^q-|\bar{v}(s,|r^-|)|^q\right)\big|}
\nonumber\\
& \leq C|r^+-|r^-||\biggl\{\int_0^1
\left(|v|^{q-2}+|\bar{v}|^{q-2}\right)|\pa_rv||v-\bar{v}|
(s,\theta r^+ + (1-\theta)|r^-|)
\,d\theta\nonumber\\
&\quad + \int_0^1|\bar{v}|^{q-1}|\pa_r(v-\bar{v})|
(s,\theta r^+ +(1-\theta)|r^-|)\,d\theta
\biggr\}
\nonumber
\end{align}
for (\ref{z27}), and so on. Here we set $r^{\pm}=r \pm (t-s)$ (see 
Sections \ref{z28}--\ref{z29} below).

\subsection{Proof of {\rm (\ref{z9})}.}
Let $r > 0$. Using the obvious pointwise estimate
\begin{equation}
 |w(s,\rho)| \leq w_1^{-1}(s,\rho) \|w_1 w\|_{L^{\infty}_{t,r}},
\nonumber
\end{equation}
we immediately get
\begin{align}
 |L|w|^{p}(t,r)|
 & \leq CLw_1^{-p}(t,r)\cdot
 \|w_1 w\|^{p}_{L^{\infty}_{t,r}}.
\nonumber
\end{align}
Thus we need to show
\begin{equation}
Lw_1^{-p} \leq Cw_2(t, r)^{-1}.
\nonumber
\end{equation} 
Define the two sets $D_1$, $D_2$ by
\begin{align}
 D_1&=\left\{(s, \rho) \,:\, 0<s<t,\ |r-(t-s)|<\rho<r+(t-s),\ s-\rho > 
 -|t-r|\right\}, \nonumber\\
 D_2&=\left\{(s, \rho) \,:\, 0<s<t,\ |r-(t-s)|<\rho<r+(t-s),\ s-\rho \leq 
 -|t-r|\right\}.\nonumber
\end{align}
Note that $D_1 = \varnothing$ if $r\geq t$. 
In what follows, we will estimate $L[\chi_{D_1}w_1^{-p}]$ and $L[\chi_{D_2}w_1^{-p}]$ separately.

\bigskip

\emph{Estimate of $L[\chi_{D_1}w_1^{-p}]$}.
We assume $t > r$, since otherwise $D_1 = \varnothing$. 
We further write $\chi_{D_1}=\chi_{D_{11}}+\chi_{D_{12}}$, where
\begin{align}
 D_{11}&=\left\{(s, \rho) \,:\, (s, \rho) \in D_1,\ s-\rho < \frac{s+\rho}{2} \right\}, \nonumber\\
 D_{12}&=\left\{(s, \rho) \,:\, (s, \rho) \in D_1,\ s-\rho \geq  \frac{s+\rho}{2} \right\}.\nonumber
\end{align}
Since $\lb \rho \rb \geq 4^{-1}\lb s+\rho \rb$ on $D_{11}$, we have
\begin{align}
 L[\chi_{D_{11}}w_1^{-p}]
 & \leq \frac{C}{r}\iint_{D_{11}}\lb s+\rho \rb^{1-p}\lb s-\rho \rb^{-\mu}\,d\rho ds,
\nonumber\\
\intertext{and applying the transformation $\tau = s+\rho$, $\sigma = s 
 - \rho$, we have}
 & \leq \frac{C}{r} \int_{t-r}^{t+r}\int_{-(t-r)}^{t-r}\lb \tau \rb^{1-p} 
 \lb \sigma \rb^{-\mu}\,d\sigma d\tau
\nonumber\\
 & \leq \frac{C}{r} \int_{t-r}^{t+r}\lb \tau \rb^{1-p}\,d\tau 
 \int_0^{t-r} \lb \sigma \rb^{-\mu}\,d\sigma.
\label{z30}
\end{align}
Meanwhile, we can check that
\begin{equation}
\label{z17}
\frac{1}{r} \int_{|t-r|}^{t+r} \lb \tau \rb^{-\kappa}\,d\tau 
 \leq C \lb t+r \rb^{-1} \lb t-r \rb^{1-\kappa} \quad \mbox{for}\quad \kappa >1.
\end{equation}
Indeed, if $r > t/2$ and $r > 1$, we easily see this inequality by $r \geq C\lb 
t+r \rb$; otherwise, we can use $t+r-|t-r| \leq 2r$ and $|t-r| \geq C(t+r)$.
Putting (\ref{z30}) and (\ref{z17}) together, we get
\begin{equation}
\label{z18}
 L[\chi_{D_{11}}w_1^{-p}] \leq Cw_2(t,r)^{-1}.
\end{equation}

We next consider $L[\chi_{D_{12}}w_1^{-p}]$. 
Making the change of variables $\tau = s+\rho$, $\sigma = s 
 - \rho$ in the integral, we obtain
\begin{align}
 L[\chi_{D_{12}}w_1^{-p}] 
&\leq \frac{C}{r}\iint_{D_{12}} 
\lb \rho \rb^{1-p}\lb s+\rho \rb^{-\mu}\,d\rho ds
\nonumber\\
 & \leq \frac{C}{r} \int_{t-r}^{t+r}\int_{\frac{t-r}{2}}^{t-r}
 \lb \tau - \sigma \rb^{1-p}\lb \tau \rb^{-\mu}\,d\sigma d\tau
\nonumber\\
 & \leq \frac{C}{r} \int_{t-r}^{t+r} \lb \tau-t+r \rb^{2-p}
\lb \tau \rb^{-\mu}\,d\tau.
\label{z19}
\end{align}

Now suppose $r > t/2$ and $r > 1$ first. Then we see
\begin{align}
 \lefteqn{\mbox{RHS of (\ref{z19})}}\nonumber\\
 & \leq \frac{C}{r}\lb t-r \rb^{-\mu}
\int_{t-r}^{2(t-r)} \lb \tau -t+r \rb^{2-p}\,d\tau% \nonumber\\
%&\quad 
+ \frac{C}{r}\int_{2(t-r)}^{t+r}\lb \tau \rb^{2-p-\mu}\,d\tau
\nonumber\\
 & \leq  C w_2(t,r)^{-1}.\nonumber
\end{align}
If, on the other hand, $r \leq t/2$ or $r \leq 1$, 
\begin{align}
\mbox{RHS of (\ref{z19})}
 &\leq 
 \frac{1}{r}\lb t-r \rb^{-\mu}\int_{t-r}^{t+r} \lb \tau-t+r \rb^{2-p}\,d\tau 
\nonumber\\
 & \leq C\frac{(1+r)^{3-p}-1}{r}\lb t+r \rb^{-\mu} \nonumber\\
 & \leq Cw_2(t,r)^{-1}.\nonumber
\end{align}
Thus we have shown $L[\chi_{D_{12}}w_1^{-p}] \leq Cw_2(t, r)^{-1}$. Together with
(\ref{z18}), we obtain
\begin{equation}
  L[\chi_{D_{1}}w_1^{-p}]  \leq Cw_2(t,r)^{-1}.
\nonumber
\end{equation}

\bigskip

\emph{Estimate of $L[\chi_{D_2}w_1^{-p}]$}.\ 
By direct computation, we have
\begin{align}
L[\chi_{D_2}w_1^{-p}]   & = 
 \frac{1}{r}\int_{|t-r|}^{t+r}\int_0^{%
 \min\{\rho-|t-r|, t+r-\rho\}}
 \lb \rho \rb^{1-p} \lb \rho-s \rb^{-\mu}\,ds d\rho
\nonumber\\
  & \leq
 \frac{C}{r}\int_{|t-r|}^{t+r} \lb \rho \rb^{2-p-\mu}\,d\rho.
\nonumber
\end{align}
Thus (\ref{z17}) gives $L[\chi_{D_2}w_1^{-p}] \leq Cw_2(t, r)^{-1}$.
\qed

\subsection{Proof of {\rm (\ref{z10})}.}\label{z28}

We set
\begin{equation}
 r^{\pm} = r \pm (t-s)
\nonumber
\end{equation}
for brevity. Then we have
\begin{align}
\label{z20}
 |K_+|v|^{q}(t,r)|
 \leq \int_0^t \left\{
\frac{r^+}{w_2(s,r^+)^{q}} + \frac{|r^-|}{w_2(s,|r^-|)^{q}}\right\}\,ds
 \cdot \|w_2 v\|_{L^{\infty}_{t,r}}^{q}.
\end{align}
In order to estimate this integral, 
we note that
\begin{align}
 & \lb t-r \rb \leq \lb s+|r^-| \rb \leq \lb s+r^+ \rb \quad \mbox{for} \quad 0 \leq s \leq t,\nonumber\\
 & \lb t-r \rb = \lb s-|r^-| \rb \quad \mbox{if} \quad r^- \geq 0.
\nonumber
\end{align}
Here the first inequality is a consequence of $\big||t-r| - s\big| \leq |t-r-s| 
= |r^-|$.

Recalling (\ref{z5}), we see that, for the first part of the integral in (\ref{z20}),
\begin{align}
\lb t-r \rb^{\mu/p}\int_0^t \frac{r^+}{w_2(s,r^+)^{q}}\,ds
&\leq \int_0^t \frac{r^+ \lb s+r^+ \rb^{\mu/p}}{w_2(s,r^+)^{q}}\,ds
\nonumber\\
&\leq \int_{r^+ <s/2} \lb r^+ \rb^{1-(p-2)q}
 \lb s+r^+ \rb^{\mu/p-\mu q}\,ds\nonumber\\
&\quad + \int_{r^+ \geq s/2} 
\lb s+r^+ \rb^{1+\mu/p-q} \lb s-r^+ \rb^{-(p-3+\mu)q}\,ds.
\nonumber
\end{align}
As for the second part, we estimate it as
\begin{align}
\lefteqn{\lb t-r \rb^{\mu/p}\int_0^t \frac{|r^-|}{w_2(s,|r^-|)^{q}}\,ds}
\nonumber\\
&\leq \int_{|r^-|<s/2} \lb r^- \rb^{1-(p-2)q}\lb s+|r^-| \rb^{\mu/p-\mu q}\,ds\nonumber\\
&\quad + \int_{r^- \leq -s/2} \lb s+|r^-| \rb^{1+\mu/p-q}\lb s-|r^-| \rb^{-(p-3+\mu)q}\,ds
\nonumber\\
&\quad + \int_{r^- \geq s/2} \frac{\lb t-r \rb^{\mu/p}}{\lb s+|r^-| \rb^{q-1}\lb s-|r^-| \rb^{(p-3+\mu)q}}\,ds.
\label{z21}
\end{align}
Since $\mu$ has chosen so that (\ref{z6}) holds, we see 
\begin{equation}
\label{z22}
 \frac{\mu}{p}+1-(p-2+\mu)q<-1.
\end{equation}
Therefore, for example,
\begin{align}
\lefteqn{\sum_{i=+,-}\int_{|r^i|<s/2} 
\lb r^i \rb^{1-(p-2)q} \lb s+|r^i| \rb^{\mu/p-\mu q}\,ds}
\nonumber\\
&\leq \sum_{i=+,-}\int_{|r^i|<s/2} \lb r^i \rb^{1-(p-2)q-\mu q+\mu/p}\,ds
\leq C. \nonumber
\end{align}
All the terms other than the last term of (\ref{z21}) can be treated similarly,
and we find that they are bounded by a constant. 
Finally, using (\ref{z22}) again, we obtain
\begin{align}
\lefteqn{\int_{r^- \geq s/2} \frac{\lb t-r \rb^{\mu/p}}{\lb s+|r^-| \rb^{q-1}\lb s-|r^-| \rb^{(p-3+\mu)q}}\,ds}\nonumber\\
& = \lb t-r \rb^{\mu/p-(p-3+\mu)q}\int_{s \geq \max\{0, 2(t-r)\}} \lb 
 r-t+2s \rb^{1-q}\,ds\nonumber\\
& \leq C\lb t-r \rb^{\mu/p-(p-3+\mu)q+2-q} \leq C.\nonumber
\end{align}
\qed

\subsection{Proof of {\rm (\ref{z11})}.}
We assume $0<r < 1$ throughout the proof of (\ref{z11}). 
To cancel out the unbounded factor $r^{-1}$, we begin with
\begin{align}
\frac{1}{r}|K_+|v|^{q}(t,r)|
 & =  \frac{1}{r}\left|\int_0^t \left\{r^+|v|^{q}(s,r^+) + r^-|v|^{q}(s,|r^-|)\right\}\,ds\right|
\nonumber \\
 & \leq \frac{1}{r}\int_0^t |r^+ + r^-||v|^{q}(s,r^+)\,ds
\nonumber\\
 &\quad + \frac{1}{r}\int_0^t |r^-|\left||v|^{q}(s,r^+)-|v|^{q}(s,|r^-|)\right|\,ds.
\label{z23}
\end{align}

Firstly, let us consider the simpler part, 
the first integral of the last expression. 
Since $r^++r^-=2r$, we have
\begin{align}
\frac{1}{r}\int_0^t |r^+ + r^-||v|^{q}(s,r^+)\,ds
  \leq 2\int_0^t w_2(s,r^+)^{-q}\,ds \cdot 
\|w_2v\|^{q}_{L^{\infty}_{t,r}}.
\nonumber
\end{align}
Using $\lb t \rb \leq \lb s+r^+ \rb$, we see that
\begin{align}
\lb t \rb^{\mu/p}\int_0^tw_2(s,r^+)^{-q}\,ds
 & \leq \int_{r^+<s/2}\lb r^+ \rb^{-(p-2)q}\lb s+r^+ \rb^{\mu/p-\mu q}\,ds
\nonumber\\
 & \quad + 
\int_{r^+ \geq s/2}\lb s-r^+ \rb^{-(p-3+\mu)}
\lb s+r^+ \rb^{\mu/p-q}\,ds\nonumber\\
 & \leq C \int_{r^+<s/2}\lb r^+ \rb^{-(p-2)q+\mu/p-\mu q}\,ds\nonumber\\
 & \quad + 
 C\int_{r^+ \geq s/2}\lb s-r^+ \rb^{-(p-3+\mu)q+\mu/p - q}\,ds.\nonumber
\end{align}
Thus we obtain the desired estimate by (\ref{z22}).

We turn to the remaining part of (\ref{z23}) next.
Noting that $|r^+-|r^-|| \leq 2r$, we obtain
\begin{align}
\lefteqn{\frac{1}{r}\int_0^t |r^-|\left||v|^{q}(s,r^+)-|v|^{q}(s,|r^-|)\right|\,ds}
\nonumber\\
&\leq C\int_0^1\int_0^t |r^-|
 |v(s, r_{\theta})|^{q-1}|\pa_r v(s, r_{\theta})|\,dsd\theta,
\label{z24}
\end{align}
where we set
\begin{equation}
r_{\theta}=\theta r^+ + (1-\theta)|r^-|.
\nonumber
\end{equation}
Now we use the weighted norms as bofore, and see that
\begin{align}
 \lefteqn{\mbox{RHS of } (\ref{z24})}
\nonumber\\
&\leq C\int_0^1\int_0^t 
\frac{|r^-|}{w_2(s,r_{\theta})^{q-1}}
\frac{1}{|r_\theta|
w_3(s,r_\theta)}\,dsd\theta
 \cdot \|w_2v\|^{q-1}_{L^{\infty}_{t,r}}
\|w_3r\pa_rv\|_{L^{\infty}_{t,r}}
\nonumber\\
&\leq C\int_0^1\int_0^t 
\frac{1}{w_2(s,r_{\theta})^{q-1}
w_3(s,r_\theta)}\,dsd\theta
 \cdot \|w_2v\|^{q-1}_{L^{\infty}_{t,r}}
\|w_3r\pa_rv\|_{L^{\infty}_{t,r}}.
\nonumber
\end{align}
In what follows, we will show that
\begin{equation}
\label{z25}
 \int_0^1\int_0^t
\frac{1}{w_2(s,r_{\theta})^{q-1}w_3(s,r_\theta)}\,dsd\theta
\leq C \lb t \rb^{-\mu/p}.
\end{equation}
We may assume $t \geq 1$, 
because otherwise (\ref{z25}) is easy (recall that we are assuming $0 < r < 1$).
Note that the region $r^- \geq 0$ is included in $t - 1 < s < t$.
Therefore,
\begin{align}
\lefteqn{\lb t \rb^{\mu/p}
\int_0^1\int_{r_- \geq 0}
\frac{1}{w_2(s,r_{\theta})^{q-1}w_3(s,r_\theta)}\,dsd\theta}
\nonumber\\
&\leq C\int_0^1\int_{r^- \geq 0} 
\frac{\lb s+r_{\theta} \rb^{\mu/p}}{\lb s+r_{\theta}\rb^{\mu (q-1)}}
\,dsd\theta\nonumber\\
&\leq C\int_0^1\int_{t-1}^t 
\,dsd\theta \leq C.\nonumber
\end{align}
When $r^- < 0$, we see that
the coefficients of $s$ in $r_{\theta}=(2\theta -1)r+t-s$ 
and $s-r_{\theta}=2s + (1-2\theta)r-t$
are independent of $\theta$. Thus we estimate the integral over 
the region $r^-<0$ as
\begin{align}
\lefteqn{\lb t \rb^{\mu/p}
\int_0^1\int_{r^-<0}
\frac{1}{w_2(s,r_{\theta})^{q-1}w_3(s,r_\theta)}\,dsd\theta}
\nonumber\\
&\leq C\int_0^1\int_{r_{\theta} < s/2,\ r^-<0} 
\lb r_{\theta} \rb^{-(p-2)(q-1)}\lb s+r_{\theta} \rb^{\mu/p-\mu (q-1)-\mu}
\,dsd\theta
\nonumber\\
&\quad + C\int_0^1\int_{s/2 \leq r_{\theta},\ r^-<0} 
\lb s-r_{\theta} \rb^{-(p-3+\mu)(q-1)-\mu}\lb s+r_{\theta} \rb^{\mu/p-(q-1)}
\,dsd\theta\nonumber\\
&\leq C\int_0^1\int_{r_{\theta} < s/2,\ r^-<0} 
\lb r_{\theta} \rb^{-(p-2)(q-1)+\mu/p-\mu (q-1)-\mu}
\,dsd\theta
\nonumber\\
&\quad + C\int_0^1\int_{s/2 \leq r_{\theta},\ r^-<0} 
\lb s-r_{\theta} \rb^{-(p-3+\mu)(q-1)-\mu+\mu/p-(q-1)}
\,dsd\theta.\nonumber
\end{align}
Since
\begin{align}
& \frac{\mu}{p}-\mu (q-1)-(p-2)(q-1) - \mu
 < p-4 < -1
\nonumber
\end{align}
by (\ref{z22}), 
the integrals above are bounded by a constant. 
We have finished the proof of (\ref{z11}).
\qed

\subsection{Proof of {\rm (\ref{z12})}}\label{z29}
We follow the same process as before.
\begin{align}
|K_-|w|^{p}(t,r)|
 \leq C\int_0^t \left\{ \frac{r^+}{w_1(s,r^+)^{p}}  
 + \frac{|r^-|}{w_1(s,|r^-|)^{p}}\right\}\,ds \cdot 
\|w_1w\|^{p}_{L^{\infty}_{t,r}}.
\nonumber
\end{align}
Both of integrands $|r^{\pm}|w_1(s,|r^{\pm}|)^{-1}$ can be treated similarly. 
We only show the estimate of the second part. 
We note that $\lb t-r \rb \leq \lb s+|r^-| \rb$ and  get
\begin{align}
\lefteqn{\lb t-r \rb^{\mu}
\int_0^t \frac{|r^-|}{w_1(s,|r^-|)^q}\,ds}
\nonumber\\
&\leq C\int_{|r^-|<s/2} \lb r^- \rb^{1-p}\,ds
 + C\int_{r^- \leq -s/2} \lb s-|r^-| \rb^{1-p}\,ds
\nonumber\\
&\quad + C\int_{r^- \geq s/2} \lb s+|r^-| \rb^{1-p}\,ds
\leq C.
\nonumber
\end{align}
\qed

%% file: uvSystem.tex
\section{Proof of Theorem \ref{b6}}

We start this section by explaining some additional notation used
in this section. 
It is necessary to define the norm for $1\leq p,\,q<\infty$
\begin{equation}
 \|v(\cdot)\|_{p,q}
  :=
  \left(
   \int_0^\infty
    \left(
    \int_{S^{n-1}}|v(r\omega)|^qdS_\omega
  \right)^{p/q}r^{n-1}dr
  \right)^{1/p}
\nonumber
\end{equation}
with an obvious modification for $p=\infty$
\begin{equation}
\|v(\cdot)\|_{\infty,q}
:=
\sup_{r>0}
 \left(
  \int_{S^{n-1}}
  |v(r\omega)|^qdS_\omega
  \right)^{1/q}
\nonumber
\end{equation}
where $r=|x|$, $\omega\in S^{n-1}$. 
These types of norms have been effectively used 
for the existence theory of solutions to fully nonlinear wave equations 
in \cite{LY}, \cite{LZ}. 
Let $N$ be a nonnegative integer and $\Psi$ a characteristic function 
of a set of ${\mathbb R}^n$. 
We define the norm
\begin{equation}
\|u(t,\cdot)\|_{\Gamma,N,p,q,\Psi}
:=\sum_{|\alpha|\leq N}
\|\Psi(\cdot)\Gamma^\alpha u(t,\cdot)\|_{p,q}.
\label{a1}
\end{equation}
For $\Psi\equiv 1$ in (\ref{a1}), we omit the subscript $\Psi$. 
If $p=q$, then we omit $q$. 
If $N=0$, then we omit both the subscripts $\Gamma$ and $N$. 
In sum, we abbreviate the notation of the norm 
$\|u(t,\cdot)\|_{\Gamma,N,p,q,\Psi}$ to 
\begin{eqnarray*}
& &
\|u(t,\cdot)\|_{\Gamma,N,p,q}, \mbox{~when~}\Psi\equiv 1,\\
& &
\|u(t,\cdot)\|_{\Gamma,N,p,\Psi}, \mbox{~when~}p=q,\\
& &
\|u(t,\cdot)\|_{\Gamma,N,p}, \mbox{~when~}p=q\mbox{~and~}\Psi\equiv 1,\\
& &
\|u(t,\cdot)\|_{p,q,\Psi}, \mbox{~when~}N=0,\\
& &
\|u(t,\cdot)\|_{p,q}, \mbox{~when~}N=0\mbox{~and~}\Psi\equiv 1,\\
& &
\|u(t,\cdot)\|_{p,\Psi}, \mbox{~when~}N=0\mbox{~and~}p=q,\\
& &
\|u(t,\cdot)\|_{p}, \mbox{~when~}N=0,\,p=q,\,\mbox{and~}\Psi\equiv 1.
\end{eqnarray*}

We find solutions to the Cauchy problem to the system of wave equations
\begin{equation}
\label{a12}
\partial_t^2 u-\Delta u=v^2,\quad
\partial_t^2 v-\Delta v=(\partial_t u)^2,\quad
t>0,\,\,x\in{\mathbb R}^n
\end{equation}
with initial data (\ref{b2})
by iteration. 
For any $T>0$, 
$\varepsilon>0$, and 
$(f,g)$, $({\tilde f}, {\tilde g})$
satisfying (\ref{b8})--(\ref{b9}), 
let us define the set of functions
\begin{align*}
Z:=&
Z(T,\varepsilon,f,g,{\tilde f},{\tilde g})\\
=&
\{\,(u,v)\in
C([0,T];H^1({\mathbb R}^n)
\times
{\dot H}^{1/4}({\mathbb R}^n)):\nonumber\\
&\qquad 
\partial_j\Gamma^\alpha u\in C([0,T];L^2({\mathbb R}^n))\quad
(|\alpha|\leq 2,\,0\leq j\leq n),\nonumber\\
&\qquad
|D|^{1/4}\Gamma^\alpha v\in C([0,T];L^2({\mathbb R}^n))\quad
(|\alpha|\leq 2),\nonumber\\
&\qquad
u(0)=\varepsilon f,\,
\partial_t u(0)=\varepsilon g,\,
v(0)=\varepsilon\tilde f,\,
\partial_t v(0)=\varepsilon\tilde g\,\},\nonumber
\end{align*}
where $|D|:=\sqrt{-\Delta}$. 
Denote by $u_0$ and $v_0$ the solutions to 
$\partial_t^2 u-\Delta u=0$ with data 
$(\varepsilon f,\varepsilon g)$ and 
$(\varepsilon\tilde f,\varepsilon\tilde g)$ 
at $t=0$, respectively. 
We obviously see from the discussion below that 
$(u_0, v_0)\in Z$. 
Therefore, the set $Z$ is nonempty. 
Let us define the sequence 
$\{(u_m,v_m)\}$ inductively by solving 
\begin{equation}
\label{a14}
\partial_t^2 u_m-\Delta u_m=(v_{m-1})^2,\quad
\partial_t^2 v_m-\Delta v_m=(\partial_t u_{m-1})^2,\quad
t>0,\,\,x\in{\mathbb R}^n
\end{equation}
with initial data 
\begin{equation}
(u_m(0),\partial_t u_m(0))
=
(\varepsilon f,\varepsilon g),
\quad
(v_m(0),\partial_t v_m(0))
=
(\varepsilon\tilde f,\varepsilon\tilde g).
\nonumber
\end{equation}
It is obvious from the discussion below that 
this is a well-defined sequence in $Z$. 

Using the quantity defined as 
$M:=\Lambda_1+\Lambda_2+\Lambda_1^2+\Lambda_2^2
+\Lambda_1\Lambda_2^2$ (see (\ref{b16}) and (\ref{b17}) for the 
definitions of $\Lambda_1$, $\Lambda_2$), we set 
$Z(2C_0M\varepsilon):=\{\,(u,v)\in Z\,;\,N((u,v))\leq 2C_0M\varepsilon\,\}$, 
where 
\begin{align*}
N((u,v))
:=&
\sup_{0<t<T}
(1+t)^{-1}
\|u(t)\|_2
+
\sum_{|\alpha|\leq 2\atop 0\leq j\leq n}
\sup_{0<t<T}
\|\partial_j\Gamma^\alpha u(t)\|_2\\
&+
\sum_{|\alpha|\leq 2}
\sup_{0<t<T}
(1+t)^{-1/12}\||D|^{1/4}\Gamma^\alpha v(t)\|_2.\nonumber
\end{align*}
For the constant $C_0$, see (\ref{a2}) below. 
The set $Z(2C_0M\varepsilon)$ is complete with respect to the metric 
$d((u,v),(u',v')):=N((u-u',v-v'))$ 
$((u,v),(u',v')\in Z(2C_0M\varepsilon))$. 

Let us first recall how to bound 
$(\Gamma^\alpha u_m,\partial_t\Gamma^\alpha u_m)|_{t=0}$ 
and 
$(\Gamma^\alpha v_m,\partial_t\Gamma^\alpha v_m)|_{t=0}$ 
for $|\alpha|\leq 2$; see, e.g., Section 4 of \cite{HWY2014}. 
We will rewrite them in terms of 
$\varepsilon f$, $\varepsilon g$, $\varepsilon\tilde f$, 
and $\varepsilon\tilde g$, and obtain 
for $0<\varepsilon<1$
\begin{align}
\lefteqn{\|u_m(0)\|_2
+\sum_{|\alpha|\leq 2\atop 0\leq j\leq n}
\|(\partial_j\Gamma^\alpha u_m)(0)\|_2}
\nonumber\\
&+
\sum_{|\alpha|\leq 2}
\bigl(
\|(|D|^{1/4}\Gamma^\alpha v_m)(0)\|_2
+
\|(|D|^{-3/4}\partial_t\Gamma^\alpha v_m)(0)\|_2
\bigr)
\leq
C_0M\varepsilon.
\label{a2}
\end{align}
We like to put off the proof of this bound 
until the end of this section, 
and see how it is useful in the iteration argument. 
The crucial point in the proof of Theorem \ref{b6}
is to prove the following.
\begin{proposition}\label{a22}
For $|\alpha|\leq 2$ and $0\leq j\leq n$, 
the following estimates hold$:$
\begin{eqnarray*}
& &
\sup_{0<t<T}
\|\partial_j\Gamma^\alpha u_m(t)\|_2\\
& &
\leq
CM\varepsilon
+
C(1+T)^{(10-3n)/6}
\biggl(
\sum_{|\beta|\leq 2}
\sup_{0<t<T}
\langle t\rangle^{-1/12}
\||D|^{1/4}\Gamma^\beta v_{m-1}(t)\|_2
\biggr)^2,\nonumber\\
& &
\sup_{0<t<T}
\langle t\rangle^{-1/12}
\||D|^{1/4}\Gamma^\alpha v_m(t)\|_2\\
& &
\leq
CM\varepsilon
+
C(1+T)^{(10-3n)/6}
\biggl(
\sum_{|\beta|\leq 2\atop 0\leq j\leq n}
\sup_{0<t<T}
\|\partial_j\Gamma^\beta u_{m-1}(t)\|_2
\biggr)^2.\nonumber
\end{eqnarray*}
Here, $C$ denotes a positive constant 
independent of $m$ and $T$.
\end{proposition}
For the proof of this proposition, we use several inequalities. 
%%%%%%%%%%%%%%%%%%%%%%%%%%%%%%%%%%%%%%%%%%%%%%%%%%%%%
\begin{proposition}\label{a25}
For any $2<q<\infty$, there exists a constant $C>0$ such that 
the inequality
\begin{equation}
\bigl\|\,
\|v(r\omega)\|_{L^2(S^{n-1})}
\,\bigr\|
_{L^q((\lambda,\infty);r^{n-1}dr)}
\leq
C
\lambda^{-(n-1)s(q)}\||D|^{s(q)}v\|_2
\nonumber
\end{equation}
holds for all $\lambda>0$, where $s(q):=1/2-1/q$.
\end{proposition}

{\it Proof}. See Theorem 2.10 of Li and Zhou \cite{LZ}. 
See also Section 4 of \cite{HWY2014} where a different proof is given.
\qed

\begin{proposition}
If $1\leq p<q<\infty$ and $1/q\geq 1/p-1/n$, 
then the inequality
\begin{equation}
\label{a4}
\|v(t,\cdot)\|_{q,\chi_1}
\leq
C(1+|t|)^{-n(1/p-1/q)}
\|v(t,\cdot)\|_{\Gamma,1,p}
\end{equation}
holds. Here $\chi_1$ denotes the characteristic function of the set 
$\{x\in{\mathbb R}^n:|x|<(1+|t|)/2\}$. 
\end{proposition}

{\it Proof}. 
See Theorem 2.9 of Li and Zhou \cite{LZ}.
\qed

\begin{proposition}\label{a23}
Suppose $n\geq 2$. $(1)$ If $1/2<s<n/2$, then 
the inequality 
\begin{equation}
\label{a5}
\sup_{r>0}
r^{(n/2)-s}
\|v(r\cdot)\|_{L^p(S^{n-1})}
\leq
C\||D|^s v\|_2
\end{equation}
holds, 
where $p$ is defined as
\begin{equation}
\frac1p=\frac12-\frac{s-\frac12}{n-1}.
\nonumber
\end{equation}
$(2)$ If $\sigma$ satisfies $1/2<1-\sigma<n/2$, 
then the solution $u$ to the inhomogeneous wave equation 
$\partial_t^2 u-\Delta u=F$ in 
${\mathbb R}^n\times (0,\infty)$ 
with data $(f,g)$ at $t=0$ satisfies
\begin{align}
\label{a6}
\||D|^\sigma u(t,\cdot)\|_2
\leq &
\||D|^\sigma f\|_2
+
\||D|^{\sigma-1}g\|_2
+
C\int_0^t\|F(\tau,\cdot)\|_{p_1,\chi_1}d\tau
\nonumber\\
&+
C\int_0^t
\langle\tau\rangle^{-(n/2)+1-\sigma}
\|F(\tau,\cdot)\|_{1,p_2,\chi_2}d\tau.
\end{align}
Here $p_1$ and $p_2$ are defined as
\begin{equation}
\frac{1}{p_1}
=\frac12+\frac{1-\sigma}{n},
\quad
\frac{1}{p_2}
=
\frac12
+\frac{\frac12-\sigma}{n-1}.
\nonumber
\end{equation}
The functions $\chi_1$ and $\chi_2$ denote 
the characteristic functions of 
$\{x\in{\mathbb R}^n:|x|<(1+\tau)/2\}$ and 
$\{x\in{\mathbb R}^n:|x|>(1+\tau)/2\}$, respectively. 
The solution $u$ also satisfies 
for $j=0,\dots,n$
\begin{equation}
\label{a7}
\|\partial_j u(t,\cdot)\|_2
\leq
\||D|f\|_2+\|g\|_2
+
\int_0^t\|F(\tau,\cdot)\|_2d\tau.
\end{equation}
\end{proposition}

{\it Proof}. 
The inequality (\ref{a5}) is an immediate consequence of 
Trace Lemma (see \cite{Ho} for $n\geq 3$ and \cite{FW} for $n\geq 2$) 
and the Sobolev embedding on the unit sphere $S^{n-1}$. 
The inequality (\ref{a6}) is also a direct consequence of 
the standard Sobolev embedding 
$\|v\|_{L^{p^*}({\mathbb R}^n)}\leq C\|v\|_{{\dot H}^s({\mathbb R}^n)}$ 
$(1/p^*=1/2-s/n)$, (\ref{a5}) and the standard duality argument. 
For details, see the proof of Theorem 2.11 of Li and Zhou \cite{LZ}. 
The inequality (\ref{a7}) is the standard estimate. 
\qed

\begin{proposition}
The following commuting relations hold$:$
\begin{eqnarray*}
& &
[\Gamma_i,\Box]=0\,\,\,\mbox{for $i=0,\dots,\nu-1$, and}\,\,\,
[L_0,\Box]=-2\Box,\\
& &
[\Gamma_j,\Gamma_k]=\sum_{l=0}^\nu C^{j,k}_l\Gamma_l,\,\,\,
j,\,k=0,\dots,\nu,\\
& &
[\Gamma_j,\partial_k]
=
\sum_{l=0}^nC^{j,k}_l\partial_l,\,\,\,j=0,\dots,\nu,\,\,k=0,\dots,n.
\end{eqnarray*}
Here $C^{j,k}_l$ denotes a constant depending on 
$j$, $k$, and $l$.
\end{proposition}

{\it Proof}. We can verify these relations by direct computations. 
\qed

\begin{remark}\label{a24}
In particular, we see by this proposition 
$$
\|\partial_t u(t,\cdot)\|_{\Gamma,2,2}
\leq
C
\sum_{0\leq j\leq n\atop|\alpha|\leq 2}
\|\partial_j \Gamma^\alpha u(t,\cdot)\|_2.
$$ 
This fact will be employed in {\rm (\ref{a8})} below.
\end{remark}

{\it Proof of Proposition {\rm \ref{a22}}.}\,\,We drop the subscript $m-1$ 
until the last step of the proof. 
We also note that, in what follows, the functions $\chi_1$ and $\chi_2$ denote 
the characteristic functions of 
$\{x\in{\mathbb R}^n:|x|<(1+\tau)/2\}$ and 
$\{x\in{\mathbb R}^n:|x|>(1+\tau)/2\}$, respectively. 

In view of (\ref{a6}) with $\sigma=1/4$ and (\ref{a7}), 
our task is to bound 
$\|v(\tau,\cdot)^2\|_{\Gamma,2,2}$, 
$\|(\partial_t u(\tau,\cdot))^2\|_{\Gamma,2,p_1,\chi_1}$, 
and 
$\|(\partial_t u(\tau,\cdot))^2\|_{\Gamma,2,1,p_2,\chi_2}$, 
where
\begin{equation}
\frac{1}{p_1}
=
\frac12
+
\frac{3}{4n},
\quad
\frac{1}{p_2}
=
\frac12
+
\frac{1}{4(n-1)}.
\nonumber
\end{equation}
Let us start with the estimate of 
$\|v(\tau,\cdot)^2\|_{\Gamma,2,2}$. 
We carry it out by 
dealing with 
$\|v(\tau,\cdot)^2\|_{\Gamma,2,2,\chi_1}$ 
and 
$\|v(\tau,\cdot)^2\|_{\Gamma,2,2,\chi_2}$, 
separately.

{\it Estimate of $\|v(\tau,\cdot)^2\|_{\Gamma,2,2,\chi_1}$}. 
Define $p^*$ and $p_3$ 
as $1/p^*=1/2-1/(4n)$ 
and 
$1/2=1/p_3+1/p^*$, respectively. 
Using the H\"older inequality, (\ref{a4}), 
and the Sobolev embedding 
${\dot H}^{1/4}({\mathbb R}^n)\hookrightarrow L^{p^*}({\mathbb R}^n)$, 
we get
\begin{eqnarray}
\lefteqn{\sum_{|\alpha|\leq |\beta|\atop |\alpha|+|\beta|\leq 2}
\|
(\Gamma^\alpha v(\tau))(\Gamma^\beta v(\tau))
\|_{2,\chi_1}
\leq
C\|v(\tau)\|_{\Gamma,1,p_3,\chi_1}
\|v(\tau)\|_{\Gamma,2,p^*}}\nonumber\\
& &
\leq
C(1+\tau)^{-n(1/p^*-1/p_3)}\|v(\tau)\|_{\Gamma,2,p^*}^2\nonumber\\
& &
\leq
C(1+\tau)^{-(n-1)/2}
\biggl(
\sum_{|\alpha|\leq 2}
\||D|^{1/4}\Gamma^\alpha v(\tau)\|_2
\bigg)^2.
\label{a10}
\end{eqnarray}
Here, we have used 
$-n(1/p^*-1/p_3)=-n(2/p^*-1/2)=-(n-1)/2$.

{\it Estimate of $\|v(\tau,\cdot)^2\|_{\Gamma,2,2,\chi_2}$}. 
We use the H\"older inequality, 
the Sobolev embedding on $S^{n-1}$, 
and Proposition \ref{a25} to get
\begin{eqnarray}
\lefteqn{\sum_{|\alpha|\leq |\beta|\atop |\alpha|+|\beta|\leq 2}
\|
(\Gamma^\alpha v(\tau))(\Gamma^\beta v(\tau))
\|_{2,\chi_2}}\nonumber\\
& &
\leq
\sum_{|\alpha|\leq 2}
\|v(\tau)\|_{4,\infty,\chi_2}
\|\Gamma^\alpha v(\tau)\|_{4,2,\chi_2}
+
\sum_{|\alpha|,|\beta|\leq 1}
\|\Gamma^\alpha v(\tau)\|_{4,\chi_2}
\|\Gamma^\beta v(\tau)\|_{4,\chi_2}\nonumber\\
& &
\leq
C\|v(\tau)\|_{\Gamma,2,4,2,\chi_2}^2
\leq
C(1+\tau)^{-(n-1)/2}
\biggl(
\sum_{|\alpha|\leq 2}
\||D|^{1/4}\Gamma^\alpha v(\tau)\|_2
\bigg)^2.
\end{eqnarray}
We have finished the estimate of $\|v(\tau)^2\|_{\Gamma,2,2}$.

{\it Estimate of $\|(\partial_t u(\tau))^2\|_{\Gamma,2,p_1,\chi_1}$}. 
Recall $1/p_1=1/2+3/(4n)$. 
Using the H\"older inequality and (\ref{a4}), we get
\begin{align}
\sum_{|\alpha|\leq 2}
\|\partial_t u(\tau)\Gamma^\alpha\partial_t u(\tau)\|_{p_1,\chi_1}
&\leq
\sum_{|\alpha|\leq 2}
\|\partial_t u(\tau)\|_{4n/3,\chi_1}
\|\Gamma^\alpha\partial_t u(\tau)\|_2
\nonumber\\
&\leq
C(1+\tau)^{-n(1/2-3/(4n))}
\|\partial_t u(\tau)\|_{\Gamma,2,2}^2.
\label{a20}
\end{align}
Moreover, we also obtain
\begin{eqnarray}
\lefteqn{\sum_{|\alpha|,|\beta|\leq 1}
\|(\Gamma^\alpha\partial_t u(\tau))
(\Gamma^\beta\partial_t u(\tau))\|_{p_1,\chi_1}}
\nonumber\\
& &
\leq
\|\partial_t u(\tau)\|_{\Gamma,1,2p_1,\chi_1}^2
\leq
C(1+\tau)^{-n/2+3/4}
\|\partial_t u(\tau)\|_{\Gamma,2,2}^2.
\end{eqnarray}
Here, we have used $-2n(1/2-1/(2p_1))=-n/2+3/4$.

{\it Estimate of $\|(\partial_t u(\tau))^2\|_{\Gamma,2,1,p_2,\chi_2}$ }. 
Recall $1/p_2=1/2+1/(4(n-1))$. 
Using the H\"older inequality and 
then the Sobolev embedding on $S^{n-1}$, we get
\begin{align}
\sum_{|\alpha|\leq 2}
\|\partial_t u(\tau)\Gamma^\alpha\partial_t u(\tau)\|_{1,p_2,\chi_2}
\leq &
\sum_{|\alpha|\leq 2}
\|\partial_t u(\tau)\|_{2,4(n-1)}
\|\Gamma^\alpha\partial_t u(\tau)\|_2 \nonumber\\
\leq &
C\|\partial_t u(\tau)\|_{\Gamma,1,2}
\|\partial_t u(\tau)\|_{\Gamma,2,2}.
\label{a9}
\end{align}
We also have by the Sobolev embedding on $S^{n-1}$
\begin{eqnarray}
& &
\sum_{|\alpha|,|\beta|\leq 1}
\|
(\Gamma^\alpha\partial_t u(\tau))(\Gamma^\beta\partial_t u(\tau))
\|_{1,p_2,\chi_2}\nonumber\\
& &
\leq
\|\partial_t u(\tau)\|_{\Gamma,1,2,2p_2}^2
\leq
C\|\partial_t u(\tau)\|_{\Gamma,2,2}^2,
\end{eqnarray}
which, together with (\ref{a9}), implies for $\sigma=1/4$
\begin{equation}
\label{a11}
\langle\tau\rangle^{-(n/2)+1-\sigma}
\|(\partial_t u(\tau))^2\|_{\Gamma,2,1,p_2,\chi_2}
\leq
C\langle\tau\rangle^{-(n/2)+3/4}
\|\partial_t u(\tau)\|_{\Gamma,2,2}^2.
\end{equation}
By (\ref{a10})--(\ref{a11}), together with Remark \ref{a24}, 
we have shown for $|\alpha|\leq 2$,
\begin{align*}
\lefteqn{\sum_{j=0}^n\|\partial_j\Gamma^\alpha u_m(t)\|_2}
\nonumber\\
&\leq 
\sum_{j=0}^n\|(\partial_j\Gamma^\alpha u_m)(0)\|_2\nonumber\\
&\quad +
C\int_0^t
\langle\tau\rangle^{-(n-1)/2+1/6}d\tau
\biggl(
\sum_{|\beta|\leq 2}
\sup_{0<t<T}
\langle t\rangle^{-1/12}
\||D|^{1/4}\Gamma^\beta v_{m-1}(t)\|_2
\biggr)^2,
\end{align*}
\begin{align}
\lefteqn{\langle t\rangle^{-1/12}\||D|^{1/4}\Gamma^\alpha v_m(t)\|_2}
\nonumber\\
&
\leq
\||D|^{1/4}((\Gamma^\alpha v_m)(0))\|_2\nonumber\\
&\quad +
C\langle t\rangle^{-1/12}
\int_0^t\langle\tau\rangle^{-n/2+3/4}d\tau
\biggl(
\sum_{|\beta|\leq 2\atop 0\leq j\leq n}
\sup_{0<t<T}
\|\partial_j\Gamma^\beta u_{m-1}(t)\|_2
\biggr)^2.
\label{a8}
\end{align}
We also note that
\begin{eqnarray*}
& &
\int_0^t
\langle\tau\rangle^{-(n-1)/2+1/6}d\tau
\leq
C(1+t)^{(10-3n)/6},\\
& &
\langle t\rangle^{-1/12}
\int_0^t\langle\tau\rangle^{-n/2+3/4}d\tau
\leq
C(1+t)^{(10-3n)/6}.
\end{eqnarray*}
Hence we have obtained
\begin{equation}
N((u_m,v_m))
\leq
C_0M\varepsilon
+
\tilde{C}(1+T)^{(10-3n)/6}
N((u_{m-1},v_{m-1}))^2
\nonumber
\end{equation}
for positive constants $C_0$, $\tilde{C}$ independent of $m$, $T$.

In the same way, we get
\begin{align*}
\lefteqn{N((u_{m+1}-u_m,v_{m+1}-v_m))}\\
&\leq
\hat{C}(1+T)^{(10-3n)/6}
\bigl(
N((u_m,v_m))+N((u_{m-1},v_{m-1}))
\bigr)\nonumber\\
&\quad 
\times
N((u_m-u_{m-1},v_m-v_{m-1}))\quad(m=1,2,\dots)\nonumber
\end{align*}
for a constant $\hat{C}>0$ independent of $m$, $T$.

If we choose $T$ and $\varepsilon$ so that 
\begin{equation}
\tilde{C}(1+T)^{(10-3n)/6}C_0M\varepsilon\leq 1,
\quad
\hat{C}(1+T)^{(10-3n)/6}(2C_0M\varepsilon)\leq\frac12
\nonumber
\end{equation}
may hold, 
then it follows from the standard argument that 
$\{(u_m,v_m)\}$ converges to the limit 
in $Z(T,\varepsilon,f,g,\tilde f,\tilde g)$, 
which means that 
for any $\varepsilon$ with 
$0<\varepsilon<\min\{1,1/(\tilde{C}C_0M),1/(4\hat{C}C_0M)\}$, 
the Cauchy problem (\ref{a12}) with initial data (\ref{b2}) admits a solution $(u,v)$ in 
$Z(T_\varepsilon,\varepsilon,f,g,\tilde f,\tilde g)$ 
satisfying 
$N((u,v))\leq 2C_0M\varepsilon$. 
Here, $T_\varepsilon$ is defined as 
$$
1+T_\varepsilon
:=
\left(
\min\left\{\frac{1}{\tilde{C}C_0M}, \frac{1}{4\hat{C}C_0M}\right\}
\right)^{6/(10-3n)}
\varepsilon^{-6/(10-3n)}.
$$
Uniqueness of solutions in 
$Z(T_\varepsilon,\varepsilon,f,g,\tilde f,\tilde g)$ 
follows from essentially the same argument. 
We have finished the proof.
\qed

\vspace{0.5cm}

{\it Proof of $(\ref{a2})$}. We follow the argument of Section 4 of \cite{HWY2014}. 
We know for $|\alpha|\leq 2$
\begin{align}
\label{a16}
&
\Gamma^\alpha u_m(0)
=
\sum_{|b|\leq 2}
\sum_{|a|\leq |b|}
C_{ab}^\alpha x^a\partial^b u_m(0),\\
&
\partial_t\Gamma^\alpha u_m(0)
=
\sum_{1\leq |b|\leq 2}
\sum_{|a|\leq |b|-1}
\tilde C_{ab}^\alpha x^a\partial^b u_m(0)
+
\sum_{|a|\leq 2 \atop |b|=2}\hat C_{ab}^\alpha x^a\partial^b\partial_t u_m(0),
\end{align}
where $x^a=x_1^{a_1}\cdots x_n^{a_n}$, 
$\partial^b=\partial_t^{b_0}\cdots\partial_n^{b_n}$; 
see (4.9)--(4.10) of \cite{HWY2014}. 
Thus we have for $|\alpha|\leq 2$
\begin{eqnarray}
\lefteqn{\sum_{j=1}^n
\|\partial_j\Gamma^{\alpha} u_m(0)\|_2
+
\|\partial_t\Gamma^\alpha u_m (0)\|_2}
\nonumber\\
&\leq&
C\sum_{j=1}^n
\sum_{|b|\leq 2}
\sum_{|a|\leq |b|}
\|\partial_j(x^a\partial^b u_m(0))\|_2\nonumber\\
& &
+C\sum_{|b|=1,2}
\sum_{|a|\leq |b|-1}
\|x^a\partial^b u_m(0)\|_2
+
C\sum_{|a|\leq 2 \atop |b|=2}
\|x^a\partial^b\partial_t u_m(0)\|_2\nonumber\\
&\leq&
C\varepsilon\Lambda_1
+
C\|\partial_t^2 u_m(0)\|_2
+
C\sum_{j=1}^n
\|x_j\partial_t^2 u_m(0)\|_2\nonumber\\
& &
+
C\sum_{|a|\leq 2}
\bigl(
\|x^a\partial_x\partial_t^2 u_m(0)\|_2
+
\|x^a\partial_t^3 u_m(0)\|_2
\bigr)\nonumber\\
&\leq&
C\varepsilon\Lambda_1
+
C\sum_{|a|\leq 1}\|x^a\Box u_m(0)\|_2\nonumber\\
& &
+
C\sum_{|a|\leq 2}
\bigl(
\|x^a\partial_x\Box u_m(0)\|_2
+
\|x^a\partial_t\Box u_m(0)\|_2
\bigr).
\label{a17}
\end{eqnarray}
Using the 1st equation of (\ref{a14}), 
we get for $|a|\leq 1$
\begin{eqnarray}
& &
\|x^a\Box u_m(0)\|_2
=
\|x^a(v_{m-1}(0))^2\|_2
=
\varepsilon^2\|x^a\tilde f^2\|_2
\nonumber\\
& &
\leq
\varepsilon^2
\|x^a\tilde f\|_{p^*}\|\tilde f\|_{4n}
\leq
C\varepsilon^2
\||D|^{1/4}(x^a\tilde f)\|_2
\||D|^{n/2-1/4}\tilde f\|_2\nonumber\\
& &
\leq
\left\{
\begin{array}{ll}
\displaystyle{C\varepsilon^2\Lambda_2
\||D|^{1/4}\tilde f\|_2^{1/2}
\||D|^{5/4}\tilde f\|_2^{1/2}},&\displaystyle{n=2},\\
\displaystyle{C\varepsilon^2\Lambda_2
\||D|^{5/4}\tilde f\|_2},&\displaystyle{n=3}
\end{array}
\right.
\nonumber\\
& &
\leq
C\varepsilon^2\Lambda_2^2,
\label{a15}
\end{eqnarray}
where $1/p^*=1/2-1/(4n)$. Moreover, we get for $j=1,\dots,n$
\begin{eqnarray*}
& &
\sum_{|a|\leq 2}
\|x^a\partial_j\Box u_m(0)\|_2
=
2\varepsilon^2
\sum_{|a|\leq 2}
\|x^a\tilde f\partial_j\tilde f\|_2\\
& &
\leq
2\varepsilon^2
\sum_{|a|\leq 1}
\|x^a\tilde f\|_{p^*}\|\partial_j\tilde f\|_{4n}
+
2\varepsilon^2
\sum_{|a|,|b|=1}
\|x^a\tilde f\|_{p^*}
\|x^b\partial_j\tilde f\|_{4n}\nonumber
\leq
C\varepsilon^2\Lambda_2^2,\nonumber
\end{eqnarray*}
where, as in (\ref{a15}), we have proceeded for $|b|=1$
\begin{equation}
\|x^b\partial_j\tilde f\|_{4n}
\leq
\left\{
\begin{array}{ll}
\displaystyle{C\||D|^{1/4}(x^b\partial_j\tilde f)\|_2^{1/2}
\||D|^{5/4}(x^b\partial_j\tilde f)\|_2^{1/2}},&\displaystyle{n=2},\\
\displaystyle{C\||D|^{5/4}(x^b\partial_j\tilde f)\|_2},&\displaystyle{n=3}
\end{array}
\right.
\nonumber
\end{equation}
and dealt with $\||D|^{5/4}(x^b\partial_j\tilde f)\|_2$ as
\begin{equation}
\||D|^{5/4}(x^b\partial_j\tilde f)\|_2
\leq
C
\biggl(
\||D|^{1/4}\partial_j\tilde f\|_2
+
\sum_{|\alpha|=2}
\||D|^{1/4}(x^b\partial_x^\alpha\tilde f)\|_2
\biggr).
\nonumber
\end{equation}
We also obtain for $|a|\leq 2$
\begin{eqnarray*}
& &
\|x^a(\partial_t\Box u_m)(0)\|_2
=
2\varepsilon^2\|x^a\tilde f\tilde g\|_2\\
& &
\leq
2\varepsilon^2
\sum_{|b|\leq 1}
\|x^b\tilde f\|_{p^*}\|\tilde g\|_{4n}
+
2\varepsilon^2
\biggl(
\sum_{|b|=1}
\|x^b\tilde f\|_{p^*}
\biggr)
\biggl(
\sum_{|b|=1}
\|x^b\tilde g\|_{4n}
\biggr)
\leq
C\varepsilon^2\Lambda_2^2,\nonumber
\end{eqnarray*}
where we have handled 
$\|x^b\tilde g\|_{4n}$ $(|b|=1)$ as
\begin{equation}
\|x^b\tilde g\|_{4n}
\leq
\left\{
\begin{array}{ll}
\displaystyle{C\||D|^{1/4}(x^b\tilde g)\|_2^{1/2}
\||D|^{5/4}(x^b\tilde g)\|_2^{1/2}},&\displaystyle{n=2},\\
\displaystyle{C\||D|^{5/4}(x^b\tilde g)\|_2},&\displaystyle{n=3}
\end{array}
\right.
\nonumber
\end{equation}
and 
$\||D|^{1/4}(x^b\tilde g)\|_2$, 
$\||D|^{5/4}(x^b\tilde g)\|_2$ as 
\begin{align*}
\||D|^{1/4}(x^b\tilde g)\|_2
&\leq
C\biggl(
\||D|^{-3/4}\tilde g\|_2
+
\sum_{|\alpha|=1}
\||D|^{-3/4}(x^b\partial_x^\alpha\tilde g)\|_2
\biggr),\\
\||D|^{5/4}(x^b\tilde g)\|_2
&\leq
C\biggl(
\||D|^{1/4}\tilde g\|_2
+
\sum_{|\alpha|=1}
\||D|^{1/4}(x^b\partial_x^\alpha\tilde g)\|_2
\biggr)\\
&\leq
C
\biggl(
\sum_{|\alpha|=1}
\||D|^{-3/4}\partial_x^\alpha\tilde g\|_2
+
\sum_{|\alpha|=2}
\||D|^{-3/4}(x^b\partial_x^\alpha\tilde g)\|_2
\biggr),\nonumber
\end{align*}
respectively. 
Summing up, we have obtained
\begin{equation}
\sum_{|\alpha|\leq 2}
\sum_{j=0}^n
\|\partial_j\Gamma^\alpha u_m(0)\|_2
\leq
C\varepsilon\Lambda_1
+
C\varepsilon^2\Lambda_2^2.
\nonumber
\end{equation}
We next prove the estimate related to $v_m(0)$. 
As in (\ref{a16})--(\ref{a17}), 
we get for $|\alpha|\leq 2$
\begin{eqnarray}
& &
\||D|^{1/4}\Gamma^\alpha v_m(0)\|_2
+
\||D|^{-3/4}\partial_t\Gamma^\alpha v_m(0)\|_2\\
& &
\leq
C\sum_{|b|\leq 2 \atop |a|\leq |b|}
\||D|^{1/4}(x^a\partial^b v_m(0))\|_2
+
C\sum_{1\leq |b|\leq 2 \atop |a|\leq |b|-1}
\||D|^{-3/4}(x^a\partial^b v_m(0))\|_2\nonumber\\
& &
+
C\sum_{|a|\leq 2 \atop |b|=2}
\||D|^{-3/4}(x^a\partial^b\partial_t v_m(0))\|_2.
\label{a18}
\end{eqnarray}
The 1st term on the right-hand side above is handled as
\begin{eqnarray*}
& &
\sum_{|b|\leq 1 \atop |a|\leq |b|}
\||D|^{1/4}(x^a\partial^b v_m(0))\|_2
+
\sum_{|b|=2 \atop |a|\leq 2}
\||D|^{1/4}(x^a\partial^b v_m(0))\|_2\\
& &
\leq
C\varepsilon\Lambda_2
+
C\sum_{|a|\leq 2}
\||D|^{1/4}(x^a\Box v_m(0))\|_2.\nonumber
\end{eqnarray*}
The 2nd term on the right-hand side of (\ref{a18}) 
is treated as
\begin{eqnarray*}
& &
\sum_{|b|=1}
\||D|^{-3/4}(\partial^b v_m(0))\|_2
+
\sum_{|b|=2 \atop |a|\leq 1}
\||D|^{-3/4}(x^a\partial^b v_m(0))\|_2\\
& &
\leq
C\varepsilon\Lambda_2
+
C\sum_{|a|\leq 1}
\||D|^{-3/4}(x^a\Box v_m(0))\|_2.\nonumber
\end{eqnarray*}
We remark that we have dealt with 
$\||D|^{-3/4}(x^a\partial_i\partial_j\tilde f)\|_2$ 
$(|a|=1)$ as 
\begin{eqnarray*}
& &
\||D|^{-3/4}\partial_i(x^a\partial_j\tilde f)\|_2
+
\||D|^{-3/4}((\partial_i x^a)\partial_j\tilde f)\|_2\\
& &
\leq
\||D|^{1/4}(x^a\partial_j\tilde f)\|_2
+
\||D|^{1/4}\tilde f\|_2
\leq
2\Lambda_2.
\end{eqnarray*}
The 3rd term on the right-hand side of (\ref{a18}) 
is treated as, for $|a|\leq 2$, 
\begin{eqnarray*}
& &
\sum_{i,j=1}^n
\varepsilon
\||D|^{-3/4}(x^a\partial_i\partial_j\tilde g)\|_2\\
& &
+
\sum_{i=1}^n
\||D|^{-3/4}(x^a\partial_i\partial_t^2 v_m(0))\|_2
+
\||D|^{-3/4}(x^a\partial_t^3 v_m(0))\|_2\nonumber\\
& &
\leq
C\varepsilon\Lambda_2
+
\sum_{i=1}^n
\||D|^{-3/4}\partial_i(x^a\partial_t^2 v_m(0))\|_2
+
\sum_{|b|\leq 1}
\||D|^{-3/4}(x^b\partial_t^2 v_m(0))\|_2\nonumber\\
& &
+
\||D|^{-3/4}(x^a\Delta\partial_t v_m(0))\|_2
+
2\||D|^{-3/4}(x^a\partial_t u_{m-1}(0)\partial_t^2 u_{m-1}(0))\|_2\nonumber\\
& &
\leq
C\varepsilon\Lambda_2
+
C\||D|^{1/4}(x^a\Box v_m(0))\|_2
+
\sum_{|b|\leq 1}
\||D|^{-3/4}(x^b\Box v_m(0))\|_2\nonumber\\
& &
+2\||D|^{-3/4}(x^a\partial_t u_{m-1}(0)\partial_t^2 u_{m-1}(0))\|_2.\nonumber
\end{eqnarray*}
Since the term $\||D|^{1/4}(x^a\Box v_m(0))\|_2$ $(|a|\leq 2)$ 
can be treated as 
\begin{eqnarray*}
& &
\sum_{j=1}^n
\||D|^{-3/4}\partial_j(x^a\Box v_m(0))\|_2\\
& &
\leq
C\sum_{|b|\leq 1}
\||D|^{-3/4}(x^b\Box v_m(0))\|_2
+
\sum_{j=1}^n
\||D|^{-3/4}(x^a\partial_j\Box v_m(0))\|_2,\nonumber
\end{eqnarray*}
we finally arrive at, for $|\alpha|\leq 2$, 
\begin{align}
\lefteqn{\||D|^{1/4}\Gamma^\alpha v_m(0)\|_2
+
\||D|^{-3/4}\partial_t\Gamma^\alpha v_m(0)\|_2}
\nonumber\\
&\leq
C\varepsilon\Lambda_2
+C\varepsilon^2
\sum_{|b|\leq 1}
\||D|^{-3/4}(x^b g^2)\|_2
+C\varepsilon^2
\sum_{|b|\leq 2 \atop 1\leq j\leq n}
\||D|^{-3/4}(x^b g\partial_j g)\|_2\nonumber\\
&\quad +
C\sum_{|b|\leq 2}
\||D|^{-3/4}(x^b\partial_t u_{m-1}(0)\partial_t^2 u_{m-1}(0))\|_2.
\label{a19}
\end{align}
The 2nd term on the right-hand side of (\ref{a19}) is treated as, 
for $|b|\leq 1$, 
\begin{eqnarray}
& &
\||D|^{-3/4}(x^b g^2)\|_2
\leq
C\|x^b g^2\|_{p_1}
\leq
C\|g\|_{4n/3}\|x^b g\|_2\nonumber\\
& &
\leq
C\|g\|_{H^1}\|x^b g\|_2
\leq C\Lambda_1^2,
\end{eqnarray}
where, as in (\ref{a20}) above, 
$p_1$ is defined as 
$1/p_1=1/2+3/(4n)$. 
Similarly, the 3rd term on the right-hand side of (\ref{a19}) 
is treated as, for $j=1,\dots,n$, 
\begin{align}
\lefteqn{\sum_{|b|\leq 1}
\||D|^{-3/4}(x^b g\partial_j g)\|_2
+\sum_{|b|=2}
\||D|^{-3/4}(x^b g\partial_j g)\|_2}
\nonumber\\
&\leq
C\sum_{|b|\leq 1}
\|g\|_{4n/3}\|x^b\partial_j g\|_2
+
C\sum_{|a|,|b|=1}
\|x^a g\|_{4n/3}\|x^b\partial_j g\|_2
\leq
C\Lambda_1^2.
\end{align}
Finally, using the 1st equation of (\ref{a14}), 
we deal with the 4th term 
on the right-hand side of (\ref{a19}) as, for $|b|\leq 2$,
\begin{align}
&
\||D|^{-3/4}(x^b\partial_t u_{m-1}(0)\Box u_{m-1}(0))\|_2
\nonumber\\
&
\quad +
\||D|^{-3/4}(x^b\partial_t u_{m-1}(0)\Delta u_{m-1}(0))\|_2\nonumber\\
&=
\varepsilon^3\||D|^{-3/4}(x^b g{\tilde f}^2)\|_2
+
\varepsilon^2\||D|^{-3/4}(x^b g\Delta f)\|_2.
\label{a26}
\end{align}
We treat the 1st term above as
\begin{eqnarray}
& &
\sum_{|b|\leq 1}
\||D|^{-3/4}(x^b g{\tilde f}^2)\|_2
+
\sum_{|b|=2}
\||D|^{-3/4}(x^b g{\tilde f}^2)\|_2\nonumber\\
& &
\leq
C
\sum_{|b|\leq 1}
\|x^b g\|_2
\|\tilde f\|_{8n/3}^2
+
C\sum_{|b|=1}
\|g\|_2\|x^b\tilde f\|_{8n/3}^2.
\end{eqnarray}
Using the Sobolev embedding 
${\dot H}^{(4n-3)/8}_2({\mathbb R}^n)\hookrightarrow L^{8n/3}({\mathbb R}^n)$, 
we can handle 
$\|x^b\tilde f\|_{8n/3}$ $(|b|\leq 1)$ as 
\begin{align}
\|x^b\tilde f\|_{8n/3}
&\leq
C\||D|^{(4n-3)/8}(x^b\tilde f)\|_2
\nonumber\\
&\leq
C\||D|^{1/4}(x^b\tilde f)\|_2^{(13-4n)/8}
\||D|^{5/4}(x^b\tilde f)\|_2^{(4n-5)/8}
\leq
C\Lambda_2.
\end{align}
We treat the 2nd term on the right-hand side of (\ref{a26}) as
\begin{equation}
\label{a21}
\sum_{|b|\leq 1}
\|x^b g\|_2
\|\Delta f\|_{4n/3}
+
\sum_{|a|,|b|=1}
\|x^a g\|_{4n/3}\|x^b\Delta f\|_2
\leq
C\Lambda_1^2.
\end{equation}
Summing up, 
we have obtained by (\ref{a19})--(\ref{a21})
\begin{eqnarray*}
& &
\sum_{|\alpha|\leq 2}
\biggl(
\||D|^{1/4}\Gamma^\alpha v_m(0)\|_2
+
\||D|^{-3/4}\partial_t\Gamma^\alpha v_m(0)\|_2
\biggr)\\
& &
\leq
C\varepsilon\Lambda_2
+
C\varepsilon^2\Lambda_1^2
+
C\varepsilon^3\Lambda_1\Lambda_2^2.\nonumber
\end{eqnarray*}
We have finished the proof of (\ref{a2}).
\qed

%% file: 0511uvSystemArx.bbl
\providecommand{\bysame}{\leavevmode\hbox to3em{\hrulefill}\thinspace}
\providecommand{\MR}{\relax\ifhmode\unskip\space\fi MR }
\providecommand{\MRhref}[2]{%
  \href{http://www.ams.org/mathscinet-getitem?mr=#1}{#2}
}
\providecommand{\href}[2]{#2}

%% file: 0511uvSystemArx.bbl
\begin{thebibliography}{10}

\bibitem{A}
R.~Agemi, \emph{Blow-up of solutions to nonlinear wave equations in two space
  dimensions}, Manuscripta Math. \textbf{73} (1991), 153--162.

\bibitem{DGM}
D.~Del~Santo, V.~Georgiev, and E.~Mitidieri, \emph{Global existence of the
  solutions and formation of singularities for a class of hyperbolic systems},
  Geometric Optics and Related Topics (F.~Colombini and N.~Lerner, eds.),
  Progress in Nonlinear Differential Equations and Their Applications, vol.~32,
  1997, pp.~117--139.

\bibitem{DM}
D.~Del~Santo and E.~Mitidieri, \emph{Blow-up of solutions of a hyperbolic
  system: the critical case}, Differential Equations \textbf{34} (1998),
  1157--1163.

\bibitem{FW}
D.~Fang and C.~Wang, \emph{Weighted {Strichartz} estimates with angular
  regularity and their applications}, Forum Math. \textbf{23} (2011), 181--205.

\bibitem{GLS}
V.~Georgiev, H.~Lindblad, and C.~D. Sogge, \emph{Weighted {Strichartz}
  estimates and global existence for semilinear wave equations}, Amer. J. Math.
  \textbf{119} (1997), 1291--1319.

\bibitem{GTZ}
V.~Georgiev, H.~Takamura, and Y.~Zhou, \emph{The lifespan of solutions to
  nonlinear systems of a high-dimensional wave equation}, Nonlinear Anal.
  \textbf{64} (2006), 2215--2250.

\bibitem{Gl2}
R.~T. Glassey, \emph{Existence in the large for $\square u = f(u)$ in two space
  dimensions}, Math. Z. \textbf{178} (1981), 233--261.

\bibitem{Gl1}
\bysame, \emph{Finite-time blow-up for solutions of nonlinear wave equations},
  Math. Z. \textbf{177} (1981), 323--340.

\bibitem{HZ}
W.~Han and Y.~Zhou, \emph{Blow up for some semilinear wave equations in
  multi-space dimensions}, Comm. Partial Differential Equations \textbf{39}
  (2014), 651--665.

\bibitem{H98}
K.~Hidano, \emph{Small data scattering for wave equations with supercritical
  nonlinearity}, Proceedings of the 23rd Sapporo Symposium on Partial
  Differential Equations (Yoshikazu Giga, ed.), 1998, available at
  http://eprints3.math.sci.hokudai.ac.jp/1232/1/53.pdf, pp.~23--30.

\bibitem{HT}
K.~Hidano and K.~Tsutaya, \emph{Global existence and asymptotic behavior of
  solutions for nonlinear wave equations}, Indiana Univ. Math. J. \textbf{44}
  (1995), 1273--1305.

\bibitem{HWY2014}
K.~Hidano, C.~Wang, and K.~Yokoyama, \emph{Combined effects of two
  nonlinearities in lifespan of small solutions to semi-linear wave equations},
  arXiv:1407.6750 [math.AP], 2014.

\bibitem{Ho}
T.~Hoshiro, \emph{On weighted ${L}^2$ estimates of solutions to wave
  equations}, J. Anal. Math. \textbf{72} (1997), 127--140.

\bibitem{J79}
F.~John, \emph{Blow-up of solutions of nonlinear wave equations in three space
  dimensions}, Manuscripta Math. \textbf{28} (1979), 235--268.

\bibitem{J81}
\bysame, \emph{Blow-up for quasilinear wave equations in three space
  dimensions}, Comm. Pure Appl. Math. \textbf{34} (1981), 29--51.

\bibitem{Kl85}
S.~Klainerman, \emph{Uniform decay estimates and the {Lorentz} invariance of
  the classical wave equation}, Comm. Pure Appl. Math. \textbf{38} (1985),
  321--332.

\bibitem{Kl87}
\bysame, \emph{Remarks on the global {Sobolev} inequalities in the {Minkowski}
  space ${{\mathbb R}}^{n+1}$}, Comm. Pure Appl. Math. \textbf{40} (1987),
  111--117.

\bibitem{KK98}
H.~Kubo and K.~Kubota, \emph{Asymptotic behaviors of radially symmetric
  solutions of $\square u =|u|^p$ for super critical values $p$ in even space
  dimensions}, Japan J. Math.(N.S.) \textbf{24} (1998), 191--256.

\bibitem{KO99}
H.~Kubo and M.~Ohta, \emph{Critical blowup for systems of semilinear wave
  equations in low space dimensions}, J. Math. Anal. Appl. \textbf{240} (1999),
  340--360.

\bibitem{KTW}
Y.~Kurokawa, H.~Takamura, and K.~Wakasa, \emph{The blow-up and lifespan of
  solutions to systems of semilinear wave equation with critical exponents in
  high dimensions}, Differential Integral Equations \textbf{25} (2012),
  363--382.

\bibitem{LY}
T.~T. Li and X.~Yu, \emph{Life-span of classical solutions to fully nonlinear
  wave equations}, Comm. Partial Differential Equations \textbf{16} (1991),
  909--940.

\bibitem{LZ}
T.~T. Li and Y.~Zhou, \emph{A note on the life-span of classical solutions to
  nonlinear wave equations in four space dimensions}, Indiana Univ. Math. J.
  \textbf{44} (1995), 1207--1248.

\bibitem{LS}
H.~Lindblad and C.~D. Sogge, \emph{Long-time existence for small amplitude
  semilinear wave equations}, Amer. J. Math. \textbf{118} (1996), 1047--1135.

\bibitem{R}
M.~A. Rammaha, \emph{Finite-time blow-up for nonlinear wave equations in high
  dimensions}, Comm. Partial Differential Equations \textbf{12} (1987),
  677--700.

\bibitem{Sc85}
J.~Schaeffer, \emph{The equation $\square u = |u|^p$ for the critical value
  $p$}, Proc. Roy. Soc. Edinburgh Sect. A. \textbf{101} (1985), 31--44.

\bibitem{S84}
T.~C. Sideris, \emph{Nonexistence of global solutions to semilinear wave
  equations in high dimensions}, J. Differential Equations \textbf{52} (1984),
  378--406.

\bibitem{Tz}
N.~Tzvetkov, \emph{Existence of global solutions to nonlinear massless {Dirac}
  system and wave equation with small data}, Tsukuba J. Math. \textbf{22}
  (1998), 193--211.

\bibitem{YZ05}
B.~T. Yordanov and Q.~S. Zhang, \emph{Finite-time blowup for wave equations
  with a potential}, SIAM J. Math. Anal. \textbf{36} (2005), 1426--1433.

\bibitem{YZ06}
\bysame, \emph{Finite time blow up for critical wave equations in high
  dimensions}, J. Funct. Anal. \textbf{231} (2006), 361--374.

\bibitem{Z95}
Y.~Zhou, \emph{{Cauchy} problem for semilinear wave equations in four space
  dimensions with small initial data}, J. Partial Differential Equations
  \textbf{8} (1995), 135--144.

\bibitem{Z01}
\bysame, \emph{Blow up of solutions to the {Cauchy} problem for nonlinear wave
  equations}, Chinese Ann. Math. Ser. B \textbf{22} (2001), 275--280.

\bibitem{ZH14}
Y.~Zhou and W.~Han, \emph{Life-span of solutions to critical semilinear wave
  equations}, Comm. Partial Differential Equations \textbf{39} (2014),
  439--451.

\end{thebibliography}
